\documentclass[12pt]{amsart}

\newcommand{\ZZ}{{\mathbb Z}}
\newcommand{\R}{{\mathbb R}}
\newcommand{\C}{{\mathbb C}}

\begin{document}

\author{V.A.~Vassiliev}

\email{vva@mi.ras.ru}

\thanks{
Research supported by the Russian Fund of Basic Investigations (project
95-01-00846a) and INTAS grant (Project \# 4373)}

\title{Monodromy of complete intersections and surface potentials}

\date{Revised version was published in 1998 with an Appendix by W.!Ebeling, and
dedicated to Professor Egbert Brieskorn on the occasion of his 60th
anniversary}

\begin{abstract}
Following Newton, Ivory and Arnold, we study the Newtonian potentials of
algebraic hypersurfaces in $\R^n$. The ramification of (analytic continuations
of) these potential depends on a monodromy group, which can be considered as a
proper subgroup of the local monodromy group of a complete intersection (acting
on a {\em twisted} vanishing homology group if $n$ is odd). Studying this
monodromy group we prove, in particular, that the attraction force of a
hyperbolic layer of degree $d$ in $\R^n$ coincides with appropriate {\em
algebraic} vector-functions everywhere outside the attracting surface if $n=2$
or $d=2$, and is non-algebraic in all domains other than the hyperbolicity
domain if the surface is generic and $(d\ge 3) \& (n\ge 3) \& (n+d \ge 8)$.

Recently W.~Ebeling has removed the last restriction $d+n \ge 8$, see his
Appendix to this article.
\end{abstract}

\maketitle

\section{Introduction}

Two famous theorems of Newton assert that

a) a homogeneous spherical layer in Euclidean space does not attract bodies
inside the sphere, and

b) exterior bodies are attracted by it to the center of the sphere as by the
point-wise particle whose mass is equal to the mass of the entire sphere.

Ivory [I] extended both these theorems to the attraction of ellipsoids, and
Arnold [A 82] extended the first of them to the attraction of arbitrary
hyperbolic hypersurfaces: such a surface does not attract the particles inside
the hyperbolicity domain; see also [G 84].

In any component of the complement of the attracting surface this attraction
force coincides with a real analytic vector-function; we investigate the
ramification of this function, in particular (following one another famous
theory of Newton, see [A 87], [AV]) the question if it is algebraic or not.

We describe the monodromy group responsible for the ramification and identify
it as a subgroup of the local monodromy group of a complex complete
intersection of codimension 2 in $\C^n$. Unlike the usual local monodromy
action, this monodromy representation is reducible: e.g. the
Newton--Ivory--Arnold theorem depends on the fact that the homology class of
the set of real points of a hyperbolic surface defines an invariant element of
this action (although this element is not equal to zero: indeed, otherwise even
the potential function of the force would be zero, and not only its gradient
field, which is wrong already in the Newton's case).

Although we consider mainly the orbit of a very special cycle, formed by all
real points of a hyperbolic polynomial, all our calculations can be applied to
more general situations, e.g. when the integration cycle is an arbitrary linear
combination of real components (maybe non-compact) of an algebraic hypersurface
in $\R^n$.

In the case of odd $n$, this group acts in a vanishing homology group with
twisted coefficients (so that the corresponding kernel form $r^{2-n}ds$ of the
potential function can be integrated correctly along its elements). In \S \ 2.3
we extend the standard facts concerning vanishing homology of complete
intersections to this group, cf. [Ph 65], [G 88].

There is a (non-formal) partition of all classes of isolated singularities of
complete intersections into series with varying dimension $n$ of the ambient
space $\C^n$ (but with the constant codimension $p$ of the complete
intersection), see [E], [AGLV]; e.g., all singularities given by $p$ generic
quadrics in the spaces $\C^n$ with different $n$ and fixed $p$ form such a
series. To any such series there corresponds a series of reflection groups,
also depending on the parameter $n$; for such $n$ that $n-p$ is even, these
groups coincide with the (standard) local monodromy groups of corresponding
singularities. The homology groups described in \S \ 2.3 fill in the gap: for
$n-p$ odd, the reflection group of the natural series coincides with the
monodromy action on such a twisted homology group of the corresponding
singularity. (In the marginal case $p=1$, all the reflection groups of the
series coincide, see [GZ], [G 88].)

This is a reason why the qualitative behavior of attraction forces in the
spaces of any dimension is essentially the same, unlike the usual situation
(see e.g.  [P], [ABG], [A 87], [AV], [V 94]) when the functions given by
similar integral representations behave in very different way in the spaces of
dimensions of different parity.

For $n=2$ and arbitrary $d$, our monodromy group is finite, thus the analytic
continuation of the attraction force is finitely-valued, in particular (by the
Riemann's existence principle) algebraic, see \S \ 5.1 below. A realistic
estimate of the number of values of this continuation is given by Theorem 4. In
particular, we get a new series of examples when the attraction force coincides
with a {\em single-valued} (rational) vector-function outside the hyperbolicity
domain, see the Corollary to Theorem 4.

For $d=2$ and arbitrary $n>2,$ the monodromy group is infinite, and the orbit
of any integration cycle lies on an ellipsoidal cylinder in the vanishing
homology space. Fortunately, the integral of the attracting charge takes zero
value on the directing plane of this cylinder, thus the number of its values
along the elements of any orbit again is finite, see \S \ 5.2.

In all the other cases (when $d \ge 3$ and $n \ge 3$) it seems likely that the
monodromy group defined by the {\it generic} algebraic surface of degree $d$ in
$\R^n$ is large enough to ensure that the Newton's integral (and any other
non-zero linear form on the space of vanishing cycles) takes an infinite number
of values on the orbit of any non-invariant vector (and the unique invariant
vector is presented by the integration cycle corresponding to the hyperbolicity
domain of an hyperbolic charge). I can prove this conjecture only if the
additional restriction $d+n \ge 8$ is
satisfied\footnote{For a complete proof, removing this restriction, see the
Appendix to this article, written by W.~Ebeling}.
\medskip

Everywhere below all the homology groups $H_*(\cdot)$ are reduced modulo a
point.

\section{Vanishing homology and local monodromy of complete intersections}

Here we recall the basic facts about the local Picard--Lefschetz theory of
isolated singularities of complete intersections (see e.g. [H], [E], [AGLV])
and extend them to the case of twisted vanishing homology groups.

\subsection{Classical theory}

Let $f:(\C^n,0) \to (\C^p,0)$ be a holomorphic map, $f=(f_1, \ldots, f_p),$ and
suppose that the variety $f^{-1}(0)$ is an {\it isolated complete intersection
singularity} (ICIS) at $0$ (i.e. it is a smooth $(n-p)$-dimensional variety in
a punctured neighborhood of $0$). Suppose that the coordinates in $\C^p$ are
chosen generically, then the map $\tilde f \equiv (f_1, \ldots, f_{p-1}): \C^n
\to \C^{p-1}$ also defines an ICIS at $0$.  Let $B$ be a sufficiently small
{\em closed} disc centered at the origin in $\C^n$, and $c=(c_1, \ldots, c_p)$
a generic point very close to the origin in $\C^p$.  The corresponding
manifolds $X_f \equiv f^{-1}(c) \cap B$ and $\tilde X_{f} \equiv \tilde
f^{-1}(c_1, \ldots, c_{p-1}) \cap B$ are called the {\it Milnor fibres} of
$f$and $\tilde f$. Their homology groups are connected by the exact sequence
\begin{equation}
\cdots \to H_{n-p+1}(\tilde X_f) \to H_{n-p+1}(\tilde X_f,X_f)
\stackrel{\partial}{\longrightarrow} H_{n-p}(X_f) \to \cdots. \label{hamm}
\end{equation}
{\bf Proposition 1} (see [M], [H]). {\it The sequence (\ref{hamm}) is trivial
outside the fragment presented here. All  groups in (\ref{hamm}) are free
Abelian. Moreover, the spaces $X_f$ and $\tilde X_f$ are homotopy equivalent to
the wedges of spheres of dimensions $n-p$ and $n-p+1$ respectively.}
\medskip

The rank of $H_{n-p}(X_f)$ is called the {\it Milnor number} of the complete
intersection $f$ and is denoted by $\mu(f)$. The Milnor numbers of all {\it
quasihomogeneous} complete intersections are calculated in [GH] (in [MO] for
$p=1$); we need the following special case of this calculation. \medskip

{\bf Proposition 2.} {\it 1. The Milnor number of a homogeneous function $f:
\C^n \to \C^1$ of degree $d$ with isolated singularity at 0 is equal to
$(d-1)^n$.

2. The Milnor number of a complete intersection $f = (f_1, f_2)$ with isolated
singularity at $0$, where the functions $f_1$ and $f_2$ are  homogeneous of
degrees $a$ and $b$ respectively, is equal to $((a-1)^nb - (b-1)^na)/(a-b)$ if
$a \ne b$, and to $(a-1)^n(an-a+1)$ if $a = b$.} \medskip

The rank $\mu(f)+ \mu(\tilde f)$ of the middle group $H_{n-p+1}(\tilde
X_f,X_f)$ of (\ref{hamm}) is equal to the number of (Morse)  critical points of
the restriction of $f_p$ on $\tilde X_f$. The generators of this group are
represented by the {\it Lefschetz thimbles} defined by the (non-intersecting)
paths in $\C^1$ connecting the non-critical value $c_p$ of this restriction
with all critical values, namely, any of these thimbles is an embedded disc
swept out by the one-parametric family of {\it vanishing spheres} lying in the
varieties $f^{-1}(c_1, \ldots, c_{p-1}, \tau),$ where $\tau$ runs over the
corresponding path in $\C^1:$ when $\tau$ tends to the endpoint (i.e. to a
critical value of this restriction) the cycles of this family contract to the
corresponding critical point. These vanishing spheres in the variety $X_f$
(which corresponds to the common starting point $c_p$ of these paths) generate
the group $H_{n-p}(X_f)$, while the elements of $H_{n-p+1}(\tilde X_f)$ define
relations among them.

\subsection{Picard--Lefschetz formula for standard homology}

Let $s \subset \C^1$ be the set of all these critical values, then the group
$\pi_1(\C^1 \setminus s)$ acts naturally on all groups of (\ref{hamm}). This
action commutes with all arrows in (\ref{hamm}) and is trivial on the left-hand
group $H_{n-p+1}(\tilde X_f)$. The action on the middle and right-hand groups
$H_{n-p+1}(\tilde X_f,X_f)$, $H_{n-p}(X_f)$ is determined by the
Picard--Lefschetz formula: a class $\delta \in H_{n-p+1}(\tilde X_f,X_f)$,
being transported along a {\em simple loop} (see [Ph 67], [V 94]) $\omega_i$,
corresponding to the path connecting $c_p$ with the $i$-th critical value,
becomes
$$\delta + (-1)^{(n-p+1)(n-p+2)/2}\langle
\partial \delta, \partial \delta_i \rangle \delta_i, $$
where $\delta_i$ is the class of the thimble defined by this path, $\partial$
is the boundary operator in (\ref{hamm}), and $\langle \cdot , \cdot \rangle$
is the intersection form in $H_{n-p}(X_f)$. In particular, a similar formula
describes the monodromy action of the same loop on $H_{n-p}(X_f)$: it sends an
element $\Delta$ of this group to
\begin{equation}
\Delta + (-1)^{(n-p+1)(n-p+2)/2}\langle \Delta, \Delta_i \rangle \Delta_i,
\label{plf1}
\end{equation}
where $\Delta_i \equiv \partial \delta_i$ is the sphere vanishing along this
path.  \medskip

{\bf Proposition 3} (see e.g. [AGV]). {\it The intersection form $\langle \cdot
, \cdot \rangle$  is symmetric if $n-p$ is even and skew-symmetric if $n-p$ is
odd. The self-intersection index of any vanishing sphere is equal to $2$ if
$n-p \equiv 0 (mod \ 4)$ and to $-2$ if $n-p \equiv 2 (mod \ 4)$.} \medskip

In particular, if $n-p$ is even, then any transportation along a simple loop
$\omega_i$ acts on the group $H_{n-p}(X_f)$ (respectively, $H_{n-p+1}(\tilde
X_f,X_f)$) as the reflection in the hyperplane orthogonal to the vector
$\Delta_i$ (respectively, $\delta_i$) with respect to the intersection form in
the homology of $X_f$ (respectively, the form induced by the boundary operator
from this intersection form). The latter action is a central extension of the
former one.
\medskip

More generally, let $F$ be a $k$-parametric deformation of $f,$ i.e. a map
$F:\C^n \times \C^k \to \C^p$ such that $F(\cdot, 0) \equiv f$. For any
$\lambda \in \C^k$ lying in a sufficiently small neighborhood $D^k$ of the
origin, denote by $f_\lambda$ the map $F(\cdot, \lambda)$ and by $\tilde
f_\lambda$ the map $\C^n \to \C^{p-1}$ given by first $p-1$ coordinate
functions of $f_\lambda$. Set $X_{f,\lambda}= f^{-1}_{\lambda} \cap B$ and
$\tilde X_{f,\lambda}= \tilde f^{-1}_{\lambda} \cap B$. If $F$  is ``not very
degenerate'' then for almost all values of $\lambda$ these varieties are smooth
(with boundaries) and have the same topological type; e.g. the varieties $X_f,
\tilde X_f$ participating in (\ref{hamm}) appear in the $p$-parametric
deformation consisting of maps $f_\lambda \equiv (f_1-\lambda_1, \ldots,
f_p-\lambda_p)$ and correspond to the particular value $(\lambda_1, \ldots,
\lambda_p) = (c_1, \ldots, c_p).$ \medskip

{\bf Definition 1.} The {\it discriminant variety} $\Sigma(F)$ of $F$ is the
set of such $\lambda \in D^k$ that the topological type of the pair of
varieties  $(\tilde X_{f,\lambda}, X_{f,\lambda})$ does not coincide with that
for all neighboring $\lambda$, i.e., either the origin in $\C^{p-1}$ is a
critical value of $\tilde f_\lambda$ or the origin in $\C^1$ is a critical
value of $f_p|_{\tilde X_{f,\lambda}}$.  An exact sequence similar to
(\ref{hamm}) appears for any $\lambda \in D^k \setminus \Sigma(F)$, as well as
the monodromy action of the group $\pi_1(D^k \setminus \Sigma(F))$ on this
sequence.

Now suppose that the deformation $F$ keeps $\tilde f$ undeformed, i.e., $\tilde
f_\lambda \equiv \tilde f$ for any $\lambda$; in particular the action of this
group on the left-hand group in (\ref{hamm}) is trivial. A standard speculation
with the Zariski's theorem (see e.g. [AGV], [V]) allows us to reduce this
action to the above-considered action of the group $\pi_1(\C^1 \setminus s)$,
and thus to the Picard--Lefschetz operators.
\medskip

There is a natural map, {\it Leray tube operation}
\begin{equation}
t: H_{n-p}(X_f) \to H_{n-p+1}(\tilde X_f \setminus X_f) \label{leray}
\end{equation}
described e.g. in [Ph 67], [AGLV], [V 94]: for any cycle $\gamma$ in
$X_\lambda$ the cycle $t(\gamma)$ is swept out by the small circles in $\tilde
X_f \setminus X_f$ which are the boundaries of the fibres of the natural
fibration of the tubular neighborhood of $X_f$.

\subsection{Twisted vanishing homology of complete intersections}

Let $L_{-1}$ (respectively, $\pm \ZZ$) be the local system on $\tilde X_f
\setminus X_f$ with the fibre $\C^1$ (respectively, $ \ZZ ^1$) such that any
loop having an odd linking number with $X_f$ acts on this fibre as
multiplication by $-1$. In particular, $L_{-1} \equiv \pm \ZZ  \otimes \C$.
\medskip

Consider the obvious homomorphism
\begin{equation}
j: H_{n-p+1}(\tilde X_f \setminus X_f, L_{-1}) \to H^{lf}_{n-p+1}(\tilde X_f
\setminus X_f, L_{-1}), \label{lf}
\end{equation}
where $H^{lf}_*(\cdot)$ denotes the homology of locally finite chains.

The Lefschetz thimbles define elements also in the right-hand group of
(\ref{lf}) and in the similar group $H^{lf}_{n-p+1}(\tilde X \setminus X_f, \pm
\ZZ ):$ indeed, they are embedded discs in $\tilde X_f\setminus X_f$ with
boundary in $X_f$,  and thus their interior parts can be lifted to an arbitrary
leaf of the local system $L_{-1}$ or $\pm \ZZ $. For any such thimble $\delta_i
\in H^{lf}_{n-p+1}(\tilde X \setminus X_f, \pm \ZZ ) $ there is an element
$\kappa_i \in H_{n-p+1}(\tilde X \setminus X_f, \pm \ZZ ), $ the {\it vanishing
cycle} defined by the same path in $\C^1$, such that $j(\kappa_i) = 2
\delta_i$, see [Ph 65] and Fig. 1, where such a cycle in one-dimensional
$\tilde X$ is shown.
\bigskip

\unitlength=1.00mm \special{em:linewidth 0.4pt} \linethickness{0.4pt}
\begin{picture}(95.00,12.00)
\put(38.00,8.00){\vector(1,1){4.00}} \put(42.00,12.00){\vector(1,0){18.00}}
\put(60.00,8.00){\vector(-1,0){18.00}} \put(42.00,8.00){\vector(-1,1){4.00}}
\put(38.00,12.00){\vector(-1,0){18.00}} \put(20.00,8.00){\vector(1,0){18.00}}
\put(20.00,10.00){\oval(6.00,4.00)[l]} \put(60.00,10.00){\oval(6.00,4.00)[r]}
\put(20.00,5.00){\makebox(0,0)[cc]{$X_f$}}
\put(60.00,5.00){\makebox(0,0)[cc]{$X_f$}} \put(60.00,10.00){\circle*{1.50}}
\put(20.00,10.00){\circle*{1.50}} \put(95.00,10.00){\makebox(0,0)[cc]{Fig. 1}}
\end{picture}

 {\bf Theorem 1.} {\it a) The homomorphism (\ref{lf}) is an isomorphism,
as well as the similar homomorphism of homology groups reduced mod $\partial
\tilde X_f,$
$$ H_{n-p+1}(\tilde X_f \setminus X_f,
\partial \tilde X_f \setminus \partial X_f; L_{-1})
\to H^{lf}_{n-p+1}(\tilde X_f \setminus X_f,
\partial \tilde X_f \setminus \partial X_f; L_{-1});$$

b) the dimensions of both groups (\ref{lf}) are equal to $\nu(f) \equiv \mu(f)+
\mu(\tilde f)$, and similar homology groups in all other dimensions are
trivial;

c) the right-hand group in (\ref{lf}) is freely generated by the Lefschetz
thimbles specified by an arbitrary {\em distinguished} (see e.g. [AGV], [AGLV])
system of paths connecting the noncritical value $c_p$ of $f_p|_{\tilde X}$
with all critical values.}
\medskip

{\bf Corollary.} {\it The left-hand group in (\ref{lf}) is generated by the
{\em vanishing cycles} defined by the same paths.} \medskip

{\it Proof of the theorem.} The fact that the map (\ref{lf}) (and also its
relative version) is isomorphic is a general algebraic fact, which is true for
all local systems $L_\alpha$ with monodromy indices $\alpha \ne 1$: this
follows from the comparison of the {\it Leray spectral sequences} (see e.g.
[GrH], \S \ III.5) calculating the indicated homology groups and applied to the
identical embedding $\tilde X_f \setminus X_f \to \tilde X_f$.

The assertion of statement b) concerning the right-hand group in (\ref{lf})
follows from the similar assertion concerning the non-twisted vanishing
homology group $H_{n-p+1}(\tilde X_f, X_f; \ZZ ) $ $\equiv
H^{lf}_{n-p+1}(\tilde X_f \setminus X_f, \ZZ)$ (see Proposition 1), the fact
that $\pm \ZZ  \otimes \ZZ_2 = \ZZ \otimes \ZZ_2 = \ZZ_2$ (the constant local
system with fibre $\ZZ_2$) and from the formula of universal coefficients. The
same reasons prove that the $\ZZ_2$-torsion of the group $H^{lf}_{*}(\tilde X_f
\setminus X_f, \pm \ZZ)$ is trivial in all dimensions.

Statement c) follows now from the fact that the images of thimbles are linearly
independent already in the group $H^{lf}_{n-p+1}(\tilde X_f \setminus X_f, \pm
\ZZ) \otimes \ZZ_2$.
\medskip

Let $\Im$ be the subgroup in $H_{n-p+1}(\tilde X_f \setminus X_f, \pm \ZZ)$
generated by vanishing cycles $\kappa_i$ defined by all possible paths
(probably it coincides with entire $H_{n-p+1}(\tilde X_f \setminus X_f, \pm
\ZZ)$).
\medskip

{\bf Lemma 1.} {\it For any elements $\alpha, \beta \in \Im, $ their
intersection index is even.}
\medskip

Indeed, this index is equal to the (well-defined) intersection index of
$\alpha$ and $j(\beta),$ and $j(\beta) \in 2 H^{lf}_{n-p+1}(\tilde X_f
\setminus X_f, \pm \ZZ)$.
\medskip

Define the bilinear form $\langle \cdot , \cdot \rangle$ on $\Im$ equal to half
this intersection index.
\medskip

{\bf Proposition 4.} {\it The form $\langle \cdot , \cdot \rangle$ is symmetric
if $n-p$ is odd and is skew-symmetric if $n-p$ is even. For any basis vanishing
cycle $\kappa_i,$ $\langle \kappa_i, \kappa_i \rangle$ is equal to 2 if
$n-p\equiv 3(mod\ 4)$ and to $-2$ if $n-p\equiv 1(mod\ 4)$.}
\medskip

In the terms of this form, the monodromy action on the group $\Im$ is defined
by the same Picard--Lefschetz formula as before: the monodromy along the simple
loop $\omega_i$ takes a cycle $\kappa$ to
\begin{equation}
\kappa + (-1)^{(n-p+1)(n-p+2)/2}\langle \kappa, \kappa_i \rangle \kappa_i.
\label{pllf}
\end{equation}

\section{Surface potentials and Newton--Ivory--Arnold
theorem}

\subsection{Potential function of a surface}

Denote by $d V$ the volume differential form in ${\R}^n$, i.e. the form $dx_1
\wedge \ldots \wedge dx_n$ in the Euclidean positively oriented coordinates
$x_1, \ldots, x_n$.

Denote by $r$ the Euclidean norm in ${\R}^n$, \ $r = (x_1^2 + \cdots +
x_n^2)^{1/2}$, and by $C_n$ the area of the unit sphere in ${\R}^n$.
\medskip

{\bf Definition 2.} The {\it elementary Newton--Coulomb potential function},
or, which is the same, the {\it standard fundamental solution of the Laplace
operator} in ${\R}^n$, is the function equal to ${{1} \over {2\pi}}\hbox{ln} \,
r$ if $n=2$, and to $-r^{2-n}/((n-2)C_n)$ if $n \ge 3$. Denote this function by
$G$. \medskip

This function can be interpreted as the potential of the force of attraction by
a particle of unit mass placed at the origin, i.e., the attraction force of
this particle is equal to $\, -\hbox{grad} \ G$.

The attraction force with which a body $K$ with density distribution $P$
attracts a particle of unit mass placed at the point $x \in {\R}^n$ is equal to
minus the gradient of the corresponding {\it potential function}, whose value
at the point $x$ is equal to the integral over $K$ of the differential form
$G(x - z)P(z)d V(z)$ (if such an integral exists). \medskip

Let  $F$ be a smooth function in Euclidean space ${\R}^n$, and $M_F$ the
hypersurface $\{F = 0\}$. Suppose that $\; \hbox{grad} \ F \ne 0$ at the points
of $M_F$, so that $M_F$ is smooth.
\medskip

{\bf Definition 3.} The {\it standard charge} $\omega_F$ on the surface $M_F$
is the differential form $d V/dF$, i.e. the $(n-1)$-form such that for any
tangent frame $(l_2, \ldots, l_n)$ of $M_F$ and a transversal vector $l_1$ the
product of the values $\omega_F(l_2, \ldots, l_n)$ and $(dF, l_1)$ is equal to
the value $dV(l_1, \ldots, l_n)$. The {\it natural orientation} of the surface
$M_F$ is the orientation defined by this differential form.
\medskip

In particular, the value at a point $x \not \in M_F$  of the limit of potential
functions of homogeneous (with density $1/\epsilon$) distributions of charges
between the surfaces $F=0$ and $F=\epsilon$ is equal to the integral of the
standard charge form
\begin{equation}
G(x - z)\omega_F(z) \label{surfpot}
\end{equation}
along the naturally oriented surface $M_F$.

In a similar way, any function $P$ on the surface $M_F$ defines the charge
$P\cdot \omega_F$, which is called the {\it standard charge with density} $P$;
the potential at the point $x$ of this charge is equal to the integral of the
form
\begin{equation}
G(x-z)P(z)\omega_F(z) \label{polpot}
\end{equation}
along the naturally oriented surface $M_F$. The attraction force of this charge
is equal to minus the gradient of this potential function.

In these terms, theorems of Newton and Ivory look as follows.
\medskip

{\bf Theorem.} {\it The potential of the standard charge of the sphere
(respectively, an ellipsoid) in ${\R}^n$ given by the canonical equation (i.e.
by a polynomial $F$ of degree 2) is equal to a constant inside the sphere (the
ellipsoid), while outside it coincides (up to multiplicative constant) with the
potential function defined by any smaller ellipsoid confocal to ours.} \medskip

Arnold extended the ``interior'' part of this theorem to all {\it hyperbolic}
layers. \medskip

{\bf Definition 4.} An algebraic hypersurface $M$ of degree $d$ in $\R P^n$ is
{\it strictly hyperbolic} with respect to a point $x \in \R P^n \setminus M$ if
any real line through $x$ intersects $M$ at exactly $d$ different real points.
A polynomial $F: \R^n \to \R$ is strictly hyperbolic with respect to the point
$x \in \R^n$ if the projective closure $\bar M_F$ of the corresponding surface
$M_F$ is.
\medskip

{\bf Proposition 5} (see e.g. [ABG]). {\it If a hypersurface $M \subset \R P^n$
is strictly hyperbolic with respect to a point $x$, then  it is also strictly
hyperbolic with respect to any point in the same component of the complement of
$M$.  Any strictly hyperbolic hypersurface is smooth.} \medskip

{\bf Definition 5.} The {\it hyperbolicity domain} of a surface $M$ is the
union  of points $x$ such that $M$ is hyperbolic with respect to $x$. \medskip

{\bf Proposition 6} (see [N]). {\it The set of all hypersurfaces $M$ of given
degree $d$ in ${\R}P^n$, which are strictly hyperbolic with respect to a given
point x, is contractible (or, equivalently, the set of all polynomials of
degree $d$ defining them consists of two contractible components).} \medskip

In particular, all the strictly hyperbolic surfaces $M$ of a given degree $d$
in $\R P^n$ are situated topologically in the same way: if $d$ is even, then
$M$ is ambient (and even rigid) isotopic to the union of $[d/2]$ concentric
spheres lying in an affine chart in ${\R}P^n$; if $d$ is odd, then $M$ is
isotopic to the union of $[d/2]$ concentric spheres plus the improper
projective hyperplane. The hyperbolicity domain consists of the interior points
of the ``most interior'' spheroid. This spheroid is always convex in ${\R}P^n$,
in particular, the hyperbolicity domain in ${\R}^n$ may consist of at most two
connected components. \medskip

The hyperbolic surface $M_F$ separates the space $\R^n$ into {\it zones}: the
$k$-th zone consists of all points $x \in \R^n \setminus M_F$ such that the
minimal number of intersection points of $M_F$ with segments connecting $x$ and
points of the hyperbolicity domains is equal to $k$. In particular, the maximal
index $k$ of a zone is equal to $[d/2]+1$ if $d$ is odd and the hyperbolicity
domain in $\R^n$ consists of one component, and is equal to $[d/2]$ otherwise.
\medskip

Given a strictly hyperbolic polynomial $F$, let us fix some path-component of
its hyperbolicity domain in ${\R}^n$, and number the components of $M_F$
starting from the boundary of this component (which becomes number 1), its
neighboring component gets number 2, etc. \medskip

{\bf Definition 6.} The {\it Arnold cycle} of $F$ is the manifold $\bar M_F$,
oriented in such a way that in the restriction to its finite part $M_F$ all odd
components are taken with the natural orientation (see Definition 3), while all
even components are taken with the reversed orientations.

The {\it hyperbolic potential} (respectively, {\it hyperbolic potential with
density $P$}) of the surface $M_F$ at a point $x \in \R^n \setminus M_F$ is the
integral of the form (\ref{surfpot}) (respectively, (\ref{polpot})) along the
Arnold cycle. As usual, the {\it attraction forces} defined by these potentials
are equal to minus the gradients of the potential functions. \medskip

{\bf Lemma 2.} {\it This definition of the Arnold cycle is correct, i.e. the
orientations of different non-compact components of $M_F$ thus defined are the
restrictions of the same orientation of the corresponding components of $\bar
M_F$.} \medskip

The proof is immediate. \medskip

{\bf Theorem} (see [A 82]). {\it The hyperbolic potential of the surface $M_F$
(and moreover any hyperbolic potential with density $P$, where $P$ is a
polynomial of degree $\le d-2$) is constant inside the hyperbolicity domain.}
\medskip

(In other words, the points of the hyperbolicity domain are not attracted by
the standard charge on $M_F$ taken with sign 1 or $-1$ depending on the parity
of the number of the component on which this charge is distributed.) \medskip

The proof follows Newton's original proof: for any infinitesimally narrow cone
centred at the point $x$, whose direction is not asymptotic for the surface
$M_F$, the forces of attraction to the pieces of $M_F$ cut by the cone
annihilate one another. Indeed, let us restrict the polynomial $F$ to the line
$L$ in ${\R}^n$ through $x$ contained in this cone; then this attraction force
is equal to the solid angle of our cone multiplied by the sum of the numbers
$P(A_i)/ F'(A_i)$ over all zeros $A_i$ of the polynomial $F|_L$. The last sum
is zero because it is the sum of the residues of a rational function over all
its complex poles. \medskip

The restriction $\deg \ P \le d-2$ from the Arnold's theorem ensures that the
integration form (\ref{polpot}) is ``regular at infinity'', i.e. extends to a
holomorphic form on the projective hypersurface $\bar M_F$. Givental [G 84]
remarked that a similar statement is true for polynomial potentials of
arbitrary degree if the integration cycle $M_F$ is compact in $\R^n$: in this
case the potential function in the hyperbolicity domain coincides with a
polynomial of degree $\le \deg \ P - d + 2$.
\medskip

In other domains the potential also coincides with real analytic functions; in
the next sections we study the global behavior of these functions, in
particular their algebraicity. The ramification of these functions is defined
by the action of certain monodromy group on a certain homology group; in the
next \S \ 4 we define these objects, and in \S \ 5 we calculate this monodromy
group.

\section{Monodromy group responsible for the ramification of potentials}

\subsection{Homology groups}

For any point $x \in {\C}^n, \ x = (x_1, \ldots, x_n),$ denote by $S(x)$ the
cone in ${\C}^n$ given by the equation
\begin{equation}
(z_1 - x_1)^2 + \cdots + (z_n - x_n)^2 = 0. \label{cone}
\end{equation}
Denote by $@ \equiv @(x)$ a local system over ${\C}^n \setminus S(x)$ with
fibre ${\ZZ}$ such that the corresponding representation $\pi_1({\C}^n
\setminus S(x)) \to \ \hbox{\rm Aut}({\ZZ})$ maps the loops whose linking
numbers with $S(x)$ are odd to the multiplication by $-1$.

We specify this local system in such a way that integrals of the form $r(\cdot\
- x)dz_1 \wedge \ldots \wedge dz_n$ along the $(n-1)$-dimensional cycles with
coefficients in it are well defined. Namely, we consider the two-fold covering
over ${\C}^n \setminus S(x)$, on which this form is single-valued, and the
direct image in ${\C}^n  \setminus S(x)$ of this bundle under the obvious
projection of this covering. The trivial ${\ZZ}$-bundle over ${\C}^n \setminus
S(x)$ is naturally included in this direct image as a subbundle; the desired
local system is the quotient bundle of these two local systems. Obviously,
integrals of the form $r(\cdot\ - x)dz_1 \wedge \ldots \wedge dz_n$ (and of its
products by all single-valued functions) along the piecewise smooth $n$-chains
with coefficients in this local system are well-defined, and if these chains
are cycles, these integrals depend only on their homology classes.
\medskip

Let $F:{\C}^n \to {\C}$  be a polynomial, $W_F \subset \C^n$  the set of its
zeros, and $\bar W_F$ the projective closure of $W_F$.

For any $x \in \C^n$ we denote by  ${\mathcal H}(x)$ the group
\begin{equation}
H_{n-1}(W_F \setminus S(x), {\ZZ}) \label{abs}
\end{equation}
in the case of even $n$, and the group
\begin{equation}
H_{n-1}(W_F \setminus S(x), \ @(x)) \label{abstw}
\end{equation}
if $n$ is odd. \medskip

Similarly, denote by ${\mathcal PH}(x)$ the group
\begin{equation}
H_{n-1}(\bar W_F \setminus \bar S(x), {\ZZ}) \label{proj}
\end{equation}
in the case of even $n$, and the group
\begin{equation}
H_{n-1}(\bar W_F \setminus \bar S(x), \ @(x)) \label{projtw}
\end{equation}
in the case of odd $n$. \medskip

{\bf Definition 7.} If the polynomial $F$ is real (i.e., $F(\R^n) \subset \R$)
and strictly hyperbolic, then the Arnold cycle defines correctly an element of
the group ${\mathcal PH}(x)$ (and even of the group ${\mathcal H}(x)$ if $M_F$
is compact); these elements are called the {\it Arnold homology classes} and
are denoted by $PA(x)$ and $A(x)$ respectively.
\medskip

In the case of odd $n$, integrals of the form (\ref{polpot}) along
$(n-1)$-chains in $W_F \setminus S(x)$ with coefficients in $@(x)$ are well
defined, and the values of these integrals along the cycles depend only on
their homology classes in the group (\ref{abstw}). Moreover, if $\; \hbox{deg}
\, P \le d-2$, and hence the form (\ref{polpot}) is regular at infinity, then
it can be integrated along the chains in $\bar W_F \setminus \bar S(x)$, and
the integrals along the cycles depend only on their classes in the group
(\ref{projtw}).

In the case of even $n > 2$ the form (\ref{polpot}) is single-valued, and no
problems with the definition of similar integrals along the elements of the
group (\ref{abs}) (or even (\ref{proj}) if $\deg \, P \le d-2$) arise, and in
the exceptional case $n = 2$, when (\ref{polpot}) is logarithmic, we remember
that we are interested not in the potential, but in its first partial
derivatives with respect to the parameter $x$ (i.e. in the components of the
attraction force vector). Therefore we integrate not the form (\ref{polpot})
but its partial derivatives
$${x_i - z_i \over (x_1-z_1)^2+(x_2-z_2)^2} P(z) \omega_F, \ i
=1, 2\, ;$$ these forms are already single-valued and there is no problem in
integrating them along the elements of the group (\ref{abs}) (or (\ref{proj})
if $\ \deg \, P \le d-2$).

\subsection{Homological bundles}

For almost all $x \in {\C}^n$ the groups ${\mathcal H}(x)$ (respectively,
${\mathcal PH}(x)$) are naturally isomorphic to one another. The set of
exceptional $x$ (for which the pair $(\bar W_F \setminus \bar S(x), W_F
\setminus S(x))$ is not homeomorphic to these for all neighboring $x'$) belongs
to a proper algebraic subvariety in ${\C}^n$ consisting of three components:

a) $W_F$ itself,

b) the set of such $x$ that $S(x)$ and $W_F$ are tangent outside $x$ in
${\C}^n$, and

c) the set of such $x$ that the projective closure of $S(x)$ in ${\C}P^n$ is
``more nontransversal'' to the closure of $W_F$ at their infinitely distant
points.

For a generic $F$ the last component is empty, and the second is irreducible
provided additionally that $n \ge  3$. \smallskip

Denote this algebraic set of all exceptional $x \in {\C}^n$ by $\Sigma(F)$.
\medskip

Consider two fibre bundles over ${\C}^n \setminus \Sigma(F)$ whose fibres over
a point $x$ are the spaces $W_F \setminus S(x)$, $\bar W_F \setminus \bar
S(x)$, and associate with them the homological bundles whose fibres over the
same point are the groups ${\mathcal H}(x)$ and ${\mathcal PH}(x)$. As usual,
the Gauss--Manin connection in these bundles defines the monodromy
representations
\begin{equation}
\pi_1({\C}^n \setminus \Sigma(F)) \to \hbox{Aut} \ {\mathcal H}(x), \label{mon}
\end{equation}
\begin{equation}
\pi_1({\C}^n \setminus \Sigma(F)) \to \hbox{Aut} \ {\mathcal PH}(x).
\label{monproj}
\end{equation}
These representations obviously commute with the natural map ${\mathcal H}(x)
\to {\mathcal PH}(x).$

Let $u$ be the potential function of the polynomial charge $P\cdot \omega_F$,
i.e. the function defined for any $x$ by the integral of the form
(\ref{polpot}) along the Arnold cycle. The ramification of (the analytic
continuation of) the function $u$ depends on the monodromy action (\ref{mon})
(respectively, (\ref{monproj})) on the Arnold element in ${\mathcal H}(x)$
(respectively, in ${\mathcal PH}(x)$).
\medskip

Namely, for any multiindex $\nu \in {\ZZ}_+^n$ ($\nu \ne 0$ if $n = 2$)
consider the linear forms
\begin{equation}
N^{(\nu)} : {\mathcal H}(x) \to {\C}, \quad PN^{(\nu)} : {\mathcal PH}(x) \to
{\C}, \label{functs}
\end{equation}
whose values on the cycle $\gamma$ are equal to the integral along $\gamma$ of
the $\nu$-th partial derivative of the form (\ref{polpot}) with respect to the
parameter $x$.  \medskip

{\bf Proposition 7.} {\it For any $\nu$ ($\ne 0$ if $n=2$) and $x \in {\R}^n
\setminus \Sigma(F)$, the $\nu$-th partial derivative of the potential function
of the standard charge of the compact hyperbolic surface $M_F$ with density $P$
is finite-valued at $x$ if and only if the linear form $N^{(\nu)}$ takes
finitely many values on the orbit of the cycle $A(x)$ under the action of the
monodromy group (\ref{mon}). If $P$ is a polynomial of degree $\le d-2-|\nu|$,
then the same is true for non-compact hyperbolic surfaces if we replace $A(x)$
by $PA(x)$, $N^{(\nu)}$ by $PN^{(\nu)}$,  and the action (\ref{mon}) by
(\ref{monproj}).} \medskip

This is a tautology.

\subsection{The invariant cycle}

In this subsection we show that for any $F$ and $x \in \C^n \setminus
\Sigma(F)$ the representation (\ref{monproj}) has an invariant vector; if $F$
is a real hyperbolic polynomial and $x$ lies in its hyperbolicity domain, then
this cycle coincides with the Arnold homology class. \medskip

Denote by $PS \subset \C P^{n-1}$ the common ``infinite'' part of all cones
$\bar S(x) \subset \C P^n$ and by $@$ the local system over $\C P^{n-1}
\setminus PS$ such that any system $@(x)$ is induced from it by the obvious
projection with center $x$.
\medskip

{\bf Proposition 8.} {\it The groups
\begin{equation}
H_{n-1}({\C}P^{n-1} \setminus PS) \label{pr1}
\end{equation}
(if $n$ is even) and
\begin{equation}
H_{n-1}({\C}P^{n-1} \setminus PS, @) \label{pr2}
\end{equation}
(if $n$ is odd) are one-dimensional. The generators of all these groups are
presented by the class of the submanifold $\R P^{n-1} \subset \C P^{n-1}
\setminus PS.$}
\medskip

The proof is elementary.
\medskip

The obvious map $\Pi : \bar W_F \setminus \bar S(x) \to \C P^{n-1} \setminus
PS$ (projection from the center $x$) is a $d$-fold ramified covering of complex
(and thus oriented) manifolds. The variety $\Pi^{-1}(\R P^{n-1})$  admits thus
an orientation ($@(x)$-orientation if $n$ is odd) induced from the chosen
orientation of $\R P^{n-1}$; denote by $\Omega(x)$ the class of this variety in
the group ${\mathcal PH}(x)$.
\medskip

{\bf Proposition 9.} {\it 1. The classes $\Omega(x)$ for different $x$
constitute a section of the homology bundle over $\C^n \setminus \Sigma(F)$
with fibres ${\mathcal PH}(x)$, which is invariant under the Gauss--Manin
connection, in particular these classes are invariant under the representation
(\ref{monproj}).

2. If $F$ is a real hyperbolic polynomial and $x$ lies in its hyperbolicity
domain, then $\Omega(x)$ coincides with the Arnold homology class $PA(x)$.}
\medskip

This follows immediately from the construction.

\subsection{Reduced Arnold class}

For an arbitrary element $\gamma$ of the group ${\mathcal PH}(x)$, the
corresponding potential function $u_\gamma(x)$ can be defined as the integral
of the form (\ref{polpot}) along the cycle $\gamma$ (if this integral exists),
in particular the usual potential $u(x)$ coincides with $u_{PA(x)}(x)$.

In this subsection we for any point $x \in \R^n \setminus \Sigma(F)$ replace
the corresponding Arnold class $PA(x)$ by another class $P\tilde A(x),$ whose
potential function $u_{P\tilde A(x)}$ ramifies in exactly the same way, but
which is more convenient because (as we shall see later)

a) it is represented by a cycle lying in the ``finite'' part $W_F \setminus
S(x)$ of $\bar W_F \setminus \bar S(x)$ and thus defining an element $\tilde
A(x)$ of the (much better studied) group ${\mathcal H}(x)$, and

b) if $n$ is even, then this element $\tilde A(x)$ can be obtained by the
``Leray tube operation'' (\ref{leray}) from a certain homology class $\alpha(x)
\in H_{n-2}(W_F \cap S(x))$, so that the action (\ref{mon}) on it is reduced to
the similar action on this more standard group.
\medskip

Indeed, it follows from Proposition 9, that if the class $\gamma' \in {\mathcal
PH}(x)$ is obtained by the Gauss--Manin connection over some path in $\C^n
\setminus \Sigma(F)$ from the Arnold cycle $PA({\bf x})$, where ${\bf x}$ is a
point in the hyperbolicity domain of a compact hyperbolic surface, then the
potential function $u_{\gamma'}(x)$ is a single-valued holomorphic function in
$\C^n \setminus \Sigma(F)$. Therefore the ramification of our integrals defined
by the class $\gamma$ coincides with that defined by the class $\gamma -
\gamma'$ (if both integrals are well-defined).

For any point $x \in {\R}^n \setminus M_F$ we choose canonically some class
$\gamma'$ obtained in this way. Namely, we choose an arbitrary point ${\bf x}
\in \R^n$ in the hyperbolicity domain (if this domain has two components in
$\R^n$, then in the component closest to $x$, i.e. such that the segment
connecting $x$ and ${\bf x}$ has $\le [d/2]$ intersections with $M_f$ ). Then
connect $x$ with ${\bf x}$ by a complex line and take the path in this line
that goes from ${\bf x}$ to $x$ along the real segment and misses any point of
$W_F$ along a small arc in the {\it lower} complex half-line with respect to
this direction (i.e. the half-line into which the vector $i\cdot({\bf x}-x)$ is
directed). See Fig. 2.
\bigskip

\unitlength=1.00mm \special{em:linewidth 0.4pt} \linethickness{0.4pt}
\begin{picture}(92.00,13.00)
\put(10.00,10.00){\line(1,0){12.00}} \put(25.00,10.00){\oval(6.00,6.00)[b]}
\put(28.00,10.00){\line(1,0){14.00}} \put(45.00,10.00){\oval(6.00,6.00)[b]}
\put(48.00,10.00){\vector(1,0){13.00}} \put(61.00,7.00){\makebox(0,0)[cc]{$x$}}
\put(11.00,7.00){\makebox(0,0)[cc]{${\bf x}$}}
\put(92.00,8.00){\makebox(0,0)[cc]{Fig. 2}} \put(25.00,10.00){\circle*{1.33}}
\put(45.00,10.00){\circle*{1.33}}
\end{picture}

For any $x \in {\R}^n \setminus M_F$, denote by $PA_{\rm hyp}(x)$ the class in
${\mathcal PH}(x)$ obtained from $PA({\bf x})$ by the Gauss--Manin connection
over this path. We are interested in the monodromy of the class $PA(x) -
PA_{\rm hyp}(x)$, which will be called the {\it reduced Arnold class} and
denoted by $P\tilde A(x)$.

\subsection{Groups ${\mathcal H}(x)$ and the vanishing homology of complete
intersections.}

We shall consider especially carefully the case when the attracting surface
$W_F$ satisfies certain genericity conditions, namely, the following ones.
\medskip

We say that two holomorphic hypersurfaces in $\C^n$ are {\it simple tangent} at
their common point, if in some local holomorphic coordinates with origin at
this point one of them is given by the equality $z_n=0$, and the second by
$z_n= z_1^2+ \cdots + z_{n-1}^2.$
\medskip

{\bf Definition 8.} The polynomial $F$ (and the corresponding hypersurface
$W_F$) is $S$-{\it generic} if the projective closure $\bar W_F$ of $W_F$ is
smooth and transversal to the improper hyperplane ${\C}P^n \setminus {\C}^n$,
its ``infinite part'' $\bar W_F \setminus {\C}^n$ is transversal in the
improper hyperplane $\C P^{n-1}$ to the standard quadric $\{z_1^2 + \cdots +
z_n^2 = 0\}$, i.e. to the boundary of any cone $S(x)$, and additionally the set
of points at which $W_F$ is {\it simple} tangent to appropriate cones $S(x)$ is
dense in the set of all points of tangency of $W_F$ and these cones at their
nonsingular points.
\medskip

The transversality conditions from this definition can be reformulated as
follows:  let $\bar F$ be the principal (of degree $d$) homogeneous part of
$F$, and $r^2 \equiv z_1^2 + \cdots + z_n^2$, then the function $\bar F$ has an
isolated singularity at $0$, and also the pair of functions $(\bar F, r^2)$
defines a (homogeneous) complete intersection with an isolated singularity at
$0$.  \medskip

{\bf Theorem 2.} {\it Suppose that the algebraic surface $W_F = \{F=0\}$ in
${\C}^n$ is $S$-generic, $\hbox{\rm deg} \ F = d$. Then for a generic $x$ the
ranks of both groups (\ref{abs}), (\ref{abstw}) (in particular, of the group
${\mathcal H}(x)$) are equal to \  $(d-1)^n + (2(d-1)^n-d)/(d-2)$\ if $d
> 2$, and to $\; 2n\;$ if $d=2$.} \medskip

Indeed, the pair of functions $(F, r^2(\cdot \ -x))$ defining the manifolds
$W_F, S(x)$ is a perturbation of the complete intersection $(\bar F, r^2)$,
changing only terms of lower degree of these polynomials. Thus the pair $(W_F,
W_F \cap S(x))$ for smooth $W_F$ and nondiscriminant $x$ is homeomorphic to the
pair $(\tilde X_f, X_f)$ from (\ref{hamm}), and the local system $@(x)$ is
isomorphic to the system $\pm \ZZ$ on $\tilde X_f \setminus X_f,$ see \S \ 2.3.
For the group (\ref{abstw}) the assertion of the theorem follows now from
Theorem 1 and Propositions 1 and 2.

Denote by $\partial W_F$ the ``infinite part'' $\bar W_F \setminus {\C}^n$ of
$\bar W_F$. Then the group (\ref{abs}) is Poincar\'e--Lefschetz dual to the
group $H_{n-1}(\bar W_F, \partial W_F \cup (\bar W_F \cap \bar S(x)))$.
Consider the homological exact sequence of the triple $(\bar W_F, \partial W_F
\cup (\bar W_F \cap \bar S(x)), \partial W_F)$. By Proposition 1 and Poincar\'e
duality in the manifolds $W_F$, $W_F \cap S(x)$, the only nontrivial fragment
in this sequence is
\begin{eqnarray}
0 \to H_{n-1}(\bar W_F, \partial W_F) \to
H_{n-1}(\bar W_F, \partial W_F \cup (\bar W_F \cap \bar S(x))) \to \nonumber \\
\to H_{n-2}(\bar W_F \cap \bar S(x), \partial W_F \cap \bar S(x)) \to 0,
\label{exact}
\end{eqnarray}
and the assertion of our theorem about the group (\ref{abs}) follows from
Proposition 2. \medskip

{\bf Remark.} It is easy to see that the map
\begin{equation}
H_{n-2}(W_F \cap S(x)) \to H_{n-1}(W_F \setminus S(x)), \label{tube}
\end{equation}
conjugate with respect to Poincar\'e dualities to the third arrow in
(\ref{exact}), coincides with the Leray tube operation (\ref{leray}), in
particular in this case this operation is monomorphic.
\medskip

So we have identified the pair $(W_F, W_F \cap S(x))$ with a standard object of
the theory of singularities of complete intersections. The pair of functions
$(F, r^2(\cdot - x))$ defining this complete intersection participates in three
important families, which depend on $n, 1$ and $n+1$ parameters respectively.
Since all of them keep the first function $F$ unmoved, we describe only the
corresponding families of second components. The first family consists of all
functions $r^2(\cdot - \tilde x)$, $\tilde x \in \C^n;$ the second of all
functions $r^2(\cdot - x) - \tau,$ $\tau \in \C,$ and the third of all
functions
\begin{equation}
\label{ttt} \rho_\lambda \equiv z_1^2+ \cdots + z_n^2+ \lambda_1 z_1+ \cdots +
\lambda_n z_n + \lambda_0.
\end{equation}

Denote the parameter space of the third deformation by $T;$ the parameter
spaces $\C^n$ and $\C^1$ of the first and second families are obviously
included in it.

Define the set $\Sigma_T$ as the set of all such points $\lambda \in T$ that
the variety $W_F \cap \{\rho_\lambda=0\}$ is not smooth; the intersection of
$\Sigma_T$ with the parameter space of the first (respectively, the second)
subfamily coincides with $\Sigma(F)$ (respectively, the set $s$ of critical
values of the restriction of $r^2(\cdot - x)$ on $W_F$, see \S \ 2.2).

By the Zariski theorem, the obvious homomorphism $\pi_1(\C^1 \setminus s) \to
\pi_1(T \setminus \Sigma_T)$ is monomorphic, in particular the monodromy group
generated by the action of the latter group in ${\mathcal H}(x)$ coincides with
the standard monodromy group of the complete intersection $(\bar F, r^2)$
considered in \S \ 2.2, 2.3.
\medskip

{\bf Definition 9.} The monodromy group defined by the Gauss--Manin
representation $\pi_1(\C^1 \setminus s) \to Aut ({\mathcal H}(x))$ (or,
equivalently, $\pi_1(T \setminus \Sigma_T) \to Aut ({\mathcal H}(x))$) is
called the {\it big} monodromy group, while the similar monodromy group defined
by the natural action (\ref{mon}) is the {\em small} one.
\medskip

Below we shall see that the small monodromy group actually is a proper subgroup
of the big one. To describe it we need several more reductions and notions.

The subgroup ${\mathcal J}(x) \subset {\mathcal H}(x)$ for any $n$ is defined
as that generated by all vanishing cycles in $W_F \setminus S(x)$ defined by
all paths in $\C^1 \setminus s$ connecting $0$ with all the points of $s$, see
\S \ 2: for even $n$ it coincides with the image of the Leray tube map
(\ref{tube}), for odd $n$ it is just the group $\Im$ described in the end of \S
\ 2.3.
\medskip

On this subgroup there is a symmetric bilinear form $\langle \cdot , \cdot
\rangle$: in the case of odd $n$ it was defined before Proposition 4 (as half
the intersection index), and in the case of even $n$ it is induced by the tube
{\it mono}morphism (\ref{tube}) from the intersection index on the group
$H_{n-2}(W_F \cap S(x))$. By Propositions 3 and 4, for any vanishing cycle
$\alpha \in {\mathcal H}(x)$ $\langle \alpha, \alpha \rangle$ is equal to 2 if
$[{n+1 \over 2}]$ is odd and to $-2$ if $[{n+1 \over 2}]$ is even. \medskip

{\bf Lemma 3.} {\it For any $n$, the action of the big monodromy group on
${\mathcal H}(x)$ preserves the subgroup ${\mathcal J}(x)$ and the bilinear
form $\langle \cdot , \cdot \rangle$ on it.}
\medskip

This follows immediately from the Picard--Lefschetz formulae (\ref{plf1}),
(\ref{pllf}).

Now suppose that the polynomial $F$ is real and hyperbolic.
\medskip

{\bf Theorem 3.} {\it For any point $x$ from the $k$-th zone of $\ {\R}^n
\setminus \nobreak \Sigma(F)$, $k \le [d/2]$, the reduced Arnold class $P\tilde
A(x) = PA(x) - PA_{\rm hyp}(x)$ can be represented by a cycle with support in
$W_F \setminus S(x)$ which is homological in ${\mathcal H}(x)$ to the sum of
$k$ pairwise orthogonal vanishing cycles. In particular, its homology class
$\tilde A(x)$ belongs to the subgroup ${\mathcal J}(x)$, and its
self-intersection index $\langle \tilde A(x), \tilde A(x) \rangle$ is equal to
$2k$ if $[{n+1 \over 2}]$ is odd, and to $-2k$ if $[{n+1 \over 2}]$ is even.}
\medskip

Indeed, these vanishing cycles are constructed as follows. If the point $y \in
\R^n \setminus M_F$ is sufficiently close to a component of $M_F$, then in a
small disc $B \subset \C^n$ centered at $y$ the pair $(W_F, S(y))$ is
diffeomorphic to the pair consisting of the plane $\{x_1=1\}$ and the cone
$S(0);$ it is easy to see that both groups $H_{n-1}(B \cap W_F \setminus S(y))$
and $H_{n-1}(B \cap W_F \setminus S(y), @(y))$ are isomorphic to $\ZZ$ and
generated by vanishing cycles defined by the one-parametric family of maps $(F,
r^2(\cdot - y)-\tau),$ $\tau \in \C^n$ (in the first case this cycle is equal
to the tube around the vanishing cycle in $W_F \cap S(y)$).
\medskip

{\bf Lemma 4} (see [V 94], Lemma 2 in \S \ III.3.4). {\it If we go from the
hyperbolicity domain along a line in $\R^n$ and traverse a component of $M_F$,
then the Arnold class corresponding to the point after the traversing is equal
to the sum of this vanishing cycle and of the similar Arnold cycle for the
point before it transported by the Gauss--Manin connection over the arc of the
path from Fig. 2 connecting them.}
\medskip

In particular, the difference $PA(x)-PA_{hyp}(x)$ for $x$ from the $k$-th zone
is homologous to the sum of $k$ vanishing cycles; by construction all these
cycles lie in the finite domain $W_F \setminus S(x)$. The homology class of
this sum in ${\mathcal H}(x)$ is exactly the promised reduced Arnold class
$\tilde A(x)$, see \S \ 4.4. It remains only to prove that these cycles are
pairwise orthogonal. To do it, consider a model hyperbolic surface: the union
of $[d/2]$ concentric close spheres of radii $1, 1+\varepsilon, \ldots ,
1+([d/2]-1)\varepsilon$ (which do not intersect one another even in the complex
domain) and, if $d$ is odd, one plane distant from these spheres.

Although this surface is not $S$-generic, the above-described construction of
the cycle $\tilde A(x)$ can be accomplished for any point $x$ in the $k$-th
zone where $k \le [d/2]$ and, if $d$ is odd and $k=[d/2],$ then $x$ lies much
closer to the exterior ovaloid than to the additional plane. Then any of our
$k$ vanishing cycles lies on the complexification of its own sphere, in
particular they do not intersect one another, and our assertion is proved for
the (very degenerate) model hyperbolic surface. We can change this surface
arbitrarily weakly so that its closure $\bar W_F$ becomes $S$-generic and
transversal to $S(x)$, but the topological shape of the pair $(W_F, S(x))$ does
not change in a large ball in ${\C}^n$ containing all our $k$ vanishing cycles.
Therefore they have zero intersection indices also for a certain generic
hyperbolic polynomial. Finally, the set of nongeneric real hyperbolic
polynomials, all whose ``nongenericity'' lies in the complex domain, has
codimension at least 2 in the space of all strictly hyperbolic polynomials,
and, by Proposition 6, the space of pairs of the form $\{$a strictly hyperbolic
polynomial $F$ of degree $d$ in $\R^n$; a point $x$ of its  $k$-th zone with $k
\le [d/2]\}$ is open and path-connected; this gives our assertion also for
arbitrary generic $F$.

\section{Description of the small monodromy group and finiteness
theorems in the cases $n=2$ and $d=2$}

\subsection{The two-dimensional case}

Let $n=2$. Denote by $\eta(F)$ the number of factors $x_1^2 + x_2^2$ in the
decomposition of the principal part $\bar F$ of the polynomial $F$ into the
simplest real factors. (Of course, if $\eta(F)>0$ then $F$ is not $S$-generic.)
\medskip

{\bf Theorem 4.} {\it The attraction force of the standard charge, distributed
on a hyperbolic curve $\{F=0\}$ of degree $d$ in $\R^2$ coincides in the $k$-th
zone with the sum of two algebraic vector-functions, any of which is $\le
(^{d-\eta(F)}_{\ \ k})$-valued. The same is true for the standard charge with
polynomial density $P$ of degree $\le d-2$.

If the hyperbolic curve $\{F=0\}$ is compact and the density function $P$ is
holomorphic, then the corresponding attraction force coincides in the $k$-th
zone with the sum of two analytic finite-valued (and even algebraic if $P$ is a
polynomial) vector-functions, any of which also is  $\le (^{d-\eta(F)}_{\ \
k})$-valued.} \medskip

{\bf Corollary.} {\it If $d$ is even and $\bar F \equiv (x_1^2 + x_2^2)^{d/2}$,
then the attraction force coincides with a rational vector-function in the
``most nonhyperbolic'' $(d/2)$-th zone.} \medskip

{\bf Example.} If $d=2,$ then $\eta(F) \ne 0$ only in the Newtonian case (when
$M_F$ is a circle). In this case the attraction force is single-valued, in all
the other irreducible cases it is 4-valued in the 1-st zone. \medskip

Proof of Theorem 4. If $n=2$, then the surface $S(x)$ consists of two complex
lines through $x$, collinear to the lines $\{x_1 = \pm i\cdot x_2\}$. The
reduced Arnold class $\tilde A(x)$ corresponding to a point $x$ from the $k$-th
zone is represented by $2k$ small circles in $W_F \setminus S(x)$ around the
intersection points of these two lines with $W_F$:  $k$ circles around the
points of any line. It follows from the construction of Arnold cycles that all
these circles close to one line are oriented in  accordance with the complex
structure of the normal bundle of this line, while close to all points of the
other they are oriented clockwise.  The total number of such intersection
points in the finite domain for any line is equal to $d-\eta(F)$.  Moving the
point $x$ in ${\C}^2 \setminus \Sigma(F)$ we can only permute these $d-\eta(F)$
circles (and, if $W_F$ is smooth, all permutations can be realized). Therefore
the orbit of the monodromy group consists of $(^{d-\eta(F)}_{\ k})^2$ elements;
this implies Theorem 4. \medskip

{\bf Remark.} Already in this case we see that the small monodromy group
actually is smaller than the big one. Indeed, the standard (``big'') monodromy
group of the complete intersection $(\bar F, r^2)$ in $\R^2$ is just the
permutation group of all $2d$ points of the Milnor fibre. In particular, the
orbit of the reduced Arnold class from the $k$-th zone under this action
consists of $(^{\ \ 2d}_{k,k,2d-2k})$ points, which is much more than
$(^d_k)^2$ provided by Theorem 4 in the case $\eta(F)=0.$
\medskip

{\it Remark about Ivory's second theorem.} Given a hyperbolic surface, do there
exist other surfaces defining the same attraction force in some exterior zone?
If yes, these surfaces define the same ramification locus of the analytic
continuations of these forces. In the case of irreducible plane curves this
locus consists of $d(d-1)$ lines tangent to $W_F$ and parallel to the line $x_1
= i\cdot x_2$ plus $d(d-1)$ lines parallel to the line $x_1 = -i\cdot x_2$. If
$d=2$, the set of curves for which these ramification loci coincide consists of
all conics inscribed in a given rectangle whose sides are parallel to these two
directions. It is easy to see that this set is one-parametric and coincides
with the family of confocal conics. For larger $d$, such copotential families
do not exist or at least are exceptional, because the number $2d(d-1)$ of
conditions  that the curves of such a family should satisfy becomes much
greater than the dimension of the space of curves.

\subsection{Reduction of the kernel of the form $\langle \cdot, \cdot
\rangle$ and the case of conical sections}

Denote by $\hbox{Ker} {\mathcal J}(x)$ the kernel of the bilinear form $\langle
\cdot,\cdot \rangle$ on the group ${\mathcal J}(x)$, i.e. the set of all
$\gamma \in {\mathcal J}(x)$ such that $\langle \gamma, \alpha \rangle =0$ for
any $\alpha$. By the Picard--Lefschetz formula, this subspace is invariant
under the monodromy action, and hence this action on the quotient lattice \quad
$\tilde {\mathcal J}(x) \equiv {\mathcal J}(x)/\hbox{Ker} {\mathcal J}(x)$
\quad is well defined. \medskip

{\bf Theorem 5.} {\it If $F$ is $S$-generic, $x \in {\C}^n \setminus
\Sigma(F)$, and $P$ a polynomial of degree $p$, then any form $N^{(\nu)}$ (see
(\ref{functs})) with $|\nu| \ge p+2-d$ takes zero value on $\ \hbox{\rm Ker}
{\mathcal J}(x)$.} \medskip

{\it Proof.} Let $n$ be even, so that ${\mathcal J}(x) = t(H_{n-2}(W_F \cap
S(x_0)))$, see (\ref{tube}). By Poincar\'e duality in $W_F \cap S(x_0)$, the
condition $\gamma \in \hbox{\rm Ker} \ {\mathcal J}(x_0)$ implies that the
cycle $t^{-1}(\gamma) \in H_{n-2}(W_F \cap S(x_0))$ is homologous in the
projective closure $\bar W_F \cap \bar S(x_0) \subset {\C}P^n$ of $W_F \cap
S(x_0)$ to a cycle which lies in the improper subspace $\bar W_F \cap \bar
S(x_0) \cap ({\C}P^n \setminus {\C}^n)$. The tube around this homology provides
the homology of $\gamma$ to some cycle belonging to $ \partial W_F \setminus
\bar S(x_0) \equiv (\bar W_F \setminus \bar S(x_0)) \cap ({\C}P^n \setminus
{\C}^n)$. The last space is an $(n-2)$-dimensional Stein manifold, thus
$\gamma$ is homologous to zero in $\bar W_F \setminus \bar S(x_0)$. On the
other hand, the forms $D^{(\nu)}_x|_{x=x_0}G(x-y)P(y)\omega_F(y)$ with $|\nu|
\ge 2+p-d$ can be extended to holomorphic forms on $\bar W_F \setminus \bar
S(x_0)$, thus their integrals along $\gamma$ are equal to zero.

In the case of odd $n$, the condition $\gamma \in Ker \ {\mathcal J}(x_0)$ also
implies that $\gamma$ is homologous in $\bar W_F \setminus \bar S(x_0)$ (as a
cycle with coefficients in $@(x_0) \otimes {\C}$) to a cycle in the improper
subspace: indeed, by Poincar\'e duality this condition implies that $\gamma$
defines a trivial element of the group $H_{n-1}^{lf}(\bar W_F \setminus S(x_0),
\partial W_F \setminus S(x_0); @(x_0))$, and hence, by the relative part of
Theorem 1a), also of the group $ H_{n-1}(W_F \setminus S(x_0),
\partial W_F \setminus S(x_0); @(x_0) \otimes {\C}))$. The rest
of the proof is the same as for even $n$. \medskip

{\bf Corollary.} {\it In the conditions of Theorem 5, the linear form
$N^{(\nu)}$ induces a form on the quotient lattice $\tilde {\mathcal J}(x)$,
and the number of different values of this form on any orbit of the monodromy
action on ${\mathcal J}(x)$ coincides with similar number for the induced form
and induced monodromy action on} $\tilde {\mathcal J}(x)$. \medskip

{\bf Theorem 6.} {\it For any $n \ge 3$ the potential of the standard charge
(\ref{surfpot}) distributed on a strictly hyperbolic surface $\{F=0\}$ of
degree $2$ in $\R^n$ coincides in the 1-st zone with an algebraic function.}
\medskip

{\bf Proposition 10.} {\it If $n$ is even, $n>2,$ and $F$ is a generic quadric
in ${\C}^n$, then the pair consisting of the corresponding lattice ${\mathcal
J}(x)$ and the bilinear form $\langle \cdot, \cdot \rangle $ on it coincides
with that defined by the extended root system $\tilde D_{n+1}$. For odd $n$
this pair is a direct sum of the lattice $\tilde D_{n+1}$ and the
$(n-1)$-dimensional lattice with zero form on it.} \medskip

This fact in the case of even $n$ and non-twisted homology is proved in [E],
and the calculation for odd $n$ is essentially the same.
\medskip

{\it Proof of Theorem 6.} If $F$ is a {\it generic} quadric, then by the
Proposition 10 the lowered form $\langle \cdot, \cdot \rangle $ on the quotient
lattice $\tilde {\mathcal J}(x)$ is isomorphic to the canonical form on the
lattice $D_{n+1}$, in particular is elliptic. Hence the orbit of any class in
this lattice (in particular of the coset of the reduced Arnold class) under the
reduced monodromy action is finite, and any linear form takes finitely many
values on it.

Finally, the {\it non-generic} quadric $F$ can be approximated by a
one-parameter family $F_\tau,$ $\tau \in (0, \epsilon]$, of generic quadrics.
The analytic continuation of the potential function $u = u(F)$ is equal to the
limit of similar continuations of potentials $u(F_\tau)$. Hence the number of
leaves of $u(F)$ is majorized by the (common) number of leaves of any of the
$u(F_\tau)$. \medskip

This proof estimates the number of leaves of potential functions of quadrics by
the numbers of elements of length $ \sqrt {-2}$ in the lattice $D_{n+1}$. As we
shall see in the next subsection, this majorization is not sharp: a more
precise upper bound is  the number of integer points in the intersection of the
sphere of radius $\sqrt{-2}$ with a certain affine sublattice of corank $1$
that does not pass through the origin.

\subsection{Principal theorem on the small monodromy group}

The obvious map $\Pi : \bar W_F \setminus \bar S(x) \to \C P^{n-1} \setminus
PS$ (see \S \ 4.3) induces a homomorphism $\Pi_*$ of the group ${\mathcal
J}(x)$ to the group (\ref{pr1}) (if $n$ is even) or (\ref{pr2}) (if $n$ is
odd). Denote by ${\mathcal M}(x)$ the kernel of this homomorphism. \medskip

{\bf Theorem 7.} {\it Suppose that the polynomial $F:\C^n \to \C$ is
$S$-generic. Then for any $x \in \C^n \setminus \Sigma(F)$

a) the map $\Pi_*$ is epimorphic, in particular ${\mathcal M}(x)$ is a
sublattice of corank 1 in ${\mathcal J}(x)$;

b) ${\mathcal M}(x)$ is spanned by all vectors but one of some basis of
vanishing cycles in ${\mathcal J}(x)$;

c) the small monodromy group in ${\mathcal J}(x)$ is generated by reflections
(with respect to the form $\langle \cdot, \cdot \rangle$) in all the basis
vanishing cycles generating ${\mathcal M}(x)$;

d) the set of these basis vanishing cycles in ${\mathcal M}(x)$ is transitive
under the action of this small monodromy group;

e) the subgroup $\hbox{\rm Ker} \ {\mathcal J}(x) \subset {\mathcal J}(x)$
belongs to ${\mathcal M}(x)$.

If $F$ is a real hyperbolic polynomial of degree $d$, $x$ a point from the
$k$-th zone, $1 \le k \le [d/2],$ and $\tilde A(x)$ the corresponding reduced
Arnold cycle, then additionally

f) $\tilde A(x)$ belongs to $k$ times the generator of the quotient group
${\mathcal J}(x)/{\mathcal M}(x) \sim \ZZ $, in particular does not belong to
${\mathcal M}(x)$;

g) the linear form $\langle \tilde A(x), \cdot \rangle$ on ${\mathcal M}(x)$ is
not trivial.}
\medskip

For the proof of this theorem and next Theorem 8 see \S \ 6.
\medskip

{\bf Corollary.} {\it The orbit of any element of ${\mathcal J}(x)$ under the
small monodromy group lies in some affine hyperplane parallel to ${\mathcal
M}(x)$.}
\medskip

Indeed, this follows from Theorem 7c) and Picard--Lefschetz formula. \medskip

{\bf Definition 10.} A polynomial $P:\C^n \to \C$ is {\it very degenerate} with
respect to $W_F$ if it is equal to $0$ at all points $y \in W_F$ at which
appropriate surfaces of the form $S(x)$ are tangent to $W_F$ at their smooth
points. \medskip

{\bf Theorem 8.} {\it Suppose that $W_F$ is $S$-generic, and the polynomial $P$
is not very degenerate with respect to $W_F$. Then there exist multiindices
$\nu \in \ZZ_+^n$ with arbitrarily large $|\nu|$ such that for a generic $x \in
\C^n \setminus \Sigma(F)$, the restriction on ${\mathcal M}(x)$ of the linear
form $N^{(\nu)}$ (see (\ref{functs})) is not trivial.}

\subsection{Main conjectures}

{\bf Conjecture 1.} {\it If the hyperbolic polynomial $F$ of degree $d \ge 3$
in $\C^n,$  $n \ge 3,$ is $S$-generic, then the potential function of the
standard charge (\ref{polpot}) with not very degenerate $P$ does not coincide
with algebraic functions in the components of $\R^n \setminus \Sigma(F)$ other
than the hyperbolicity domain; moreover, the same is true for some arbitrarily
high partial derivatives of this potential function.} \medskip

Theorem 7 reduces this conjecture to the following Conjecture 2
(proved recently by W.~Ebeling, see the Appendix).
\medskip

{\bf Definition 11.} A triple $\ (A;  \langle \cdot,\cdot \rangle ; g) \ $
consisting of an integer lattice $A$, an {\it even} integer-valued symmetric
bilinear form on it and a group $g \subset \hbox{Aut}(A)$ generated by the
reflections in hyperplanes orthogonal to several elements $a_i$ of length
$\sqrt{-2}$ in $A$, is called {\it completely infinite} if for any element $a
\in A$ such that not all numbers $\langle a, a_i \rangle$ are equal to $0$, any
nonzero linear form $A \otimes {\C} \to {\C}$ takes infinitely many  values on
the orbit of $a$ under the action of the group $g$. \medskip

{\bf Conjecture 2.} {\it For any $S$-generic polynomial $F$ of degree $d \ge 3$
in ${\C}^n$, $n\ge 3$, the triple consisting of the group ${\mathcal M}(x)$,
the bilinear form equal (up to sign if $[{n+1 \over 2}]$ is odd) to the form
$\langle \cdot, \cdot \rangle$ defined before Lemma 3, and the ``small''
monodromy group on ${\mathcal M}(x)$, is completely infinite.}
\medskip

In [V 94] this conjecture was proved if additionally $n+d \ge 8$.
\medskip

{\bf Proposition 11.} {\it Conjecture 2 implies Conjecture 1.}
\medskip

{\it Proof.} Let $x$ be a nondiscriminant point in the $k$-th zone, $1 \le k
\le [d/2]$, for which the assertion of Theorem 8 with a certain $\nu$ is
satisfied. By Theorem 7b), g) there is a vanishing cycle $\Gamma \in {\mathcal
M}(x)$ such that $\langle \tilde A, \Gamma \rangle \ne 0$. By the
Picard--Lefschetz formula, the monodromy along the corresponding simple loop
takes $\tilde A$ to $\tilde A + \lambda \Gamma$, $\lambda \ne 0$. By Conjecture
2, for generic $x$ the form $N^{(\nu)}$ takes infinitely many values on the
orbit of the added term $\lambda \Gamma$ under the action of the small
monodromy group. On the other hand, this infinite number is estimated from
above by the number $q(q-1),$ where $q$ is the number of values of the form
$N^{\nu}$ on the orbit of $\tilde A(x)$, in particular this number $q$  is also
infinite.

Finally, for the points $x$ from the $([d/2]+1)$-th zone (if it exists) the
assertion of the Conjecture 1 follows from the fact that the potential function
defined by the charge (\ref{polpot}) obviously extends to an analytic function
on ${\R}P^n \setminus M_F$, hence its algebraicity in the $([d/2]+1)$-th zone
is equivalent to that in the zone separated from it by a piece of the improper
subspace in ${\R}P^n$; the number of the latter zone is surely less than
$[d/2]+1$.

\section{Proof of Theorems 7, 8}

All the main characters of statements a)--e) of Theorem 7 corresponding to all
$S$-generic $F$ of the same degree in $\C^n$ and all $x \not \in \Sigma(F)$ are
isomorphic to one another, therefore we can assume that $F$ is a real
hyperbolic polynomial and $x$ a real point.

The proof of statement e) follows immediately from that of Theorem 5.

Any induction step from the proof of Theorem 3 obviously increases the image of
$\tilde A(x)$ under the map $\Pi_*$ by a generator of the target homology
group; all such $k$ steps are locally topologically equivalent, and hence add a
fixed generator of this target group with the same sign. This proves statement
f) of Theorem 7, and statement a) is a direct corollary of it.

For any $k = 1, 2, \ldots, [d/2]$, and any point $x$ in the $k$-th zone,
consider the difference of the projective Arnold class $PA(x)$ and the element
in ${\mathcal PH}(x)$ obtained as in the definition of the reduced Arnold
cycles (i.e. by transportation along an arc in the lower complex half-line)
from a similar class $PA(x')$, $x'$ in the $(k-1)$-st zone. By Lemma 4, if $x$
and $x'$ are sufficiently close to one another and to the $k$-th component of
$M_F$ separating them, then this class can be realized by a cycle lying in a
small disc $B$ containing both these points $x, x'$. Denote by $a(x)$ the class
of this cycle in the group ${\mathcal H}(x)$; by continuity this class $a(x)$
is well defined also for arbitrary $x$ from the same zone (not necessarily
close to $M_F$). By Lemma 4, for all $x$ not in the hyperbolicity domain the
corresponding maps $\Pi_*$ send the elements $a(x)$ into the same element of
the group (\ref{pr1}) or (\ref{pr2}). \medskip

{\bf Theorem 9.} {\it If $n>2$, then

a) all classes $a(x)$, corresponding to all points $x \in {\R}^n \setminus
\Sigma(F)$, $x$ not in the hyperbolicity domain or in the $([d/2]+1)$-th zone,
can be obtained from one another by the Gauss--Manin connection in the homology
bundle $\{ {\mathcal H}(x) \to x \}$ over some path in ${\C}^n \setminus
\Sigma(F)$. These classes $a(x)$ do not belong to ${\mathcal M}(x)$, and any of
them, being added to the set of $\; \hbox{\rm dim} \, {\mathcal J}(x)-1$ basis
elements of ${\mathcal M}(x)$, mentioned in statements b), c) of Theorem 7,
completes this set to a basis in ${\mathcal J}(x)$;
\smallskip

b) for arbitrary $x$ in the $k$-th zone, $1 \le k \le [d/2]$, the linear form
$\langle a(x), \cdot \rangle$ on the group ${\mathcal M}(x)$, defined by our
bilinear form, is nontrivial.}

\subsection{Comparison of big and small monodromy groups}

Now we compare the fundamental groups of ${\C}^n \setminus \Sigma(F)$ and of
the complement of the discriminant variety $\Sigma_T$ of the deformation
(\ref{ttt}) of the complete intersection $(\bar F, r^2)$. Since $F$ is
$S$-generic, the set $\Sigma(F)$ consists of only two components, $W_F$ and the
set of $x \not \in W_F$ such that $S(x)$ is tangent to $W_F$; if $n>2$, then
the latter component is irreducible.

Let us choose the distinguished point ${\bf x}$ of the space $T \setminus
\Sigma_T$ in the hyperbolicity domain of the subspace ${\R}^n \setminus
\Sigma(F)$ . The group $\pi_1(T \setminus \Sigma_T)$ acts in the usual way on
the group ${\mathcal H}({\bf x})$ and generates the ``big'' monodromy group,
see \S \ 4.5.

Let $\Lambda$ be a generic 2-plane in $T$, and $L = \Lambda \cap {\C}^n$; $\bar
U$ a small neighbourhood of $L$ in the projective compactification of $T$, and
$U = \bar U \cap T$ the affine part of $\bar U$. Let $L'$ be a generic line in
$\Lambda$ through ${\bf x}$ sufficiently close to $L$, so that $L' \subset U$
and $L'$ intersects $\Sigma_T$ transversally. \medskip

{\bf Lemma 5.} {\it The obvious maps $\pi_1(L\setminus \Sigma(F)) \to
\pi_1({\C}^n \setminus \Sigma(F))$ and \\ $\pi_1(L' \setminus \Sigma_T) \to
\pi_1(U \setminus \Sigma_T) \to \pi_1(T \setminus \Sigma_T)$ are epimorphic.}
\medskip

The proof follows directly from the generalized Lefschetz theorem (see [GM]).
\medskip

Thus the small and big monodromy groups are generated by simple loops lying in
$L\setminus \Sigma_T$ and $L' \setminus \Sigma_T$, respectively. Let us compare
these collections of loops. \medskip

{\bf Lemma 6.} {\it The group ${\mathcal J}({\bf x})$ is generated by the
cycles vanishing along the paths of an arbitrary distinguished system in $L'$
connecting the distinguished point ${\bf x}$ with all points of $L' \cap
\Sigma_T$.} \medskip

Indeed, the group $\pi_1(L' \setminus \Sigma_T)$ acts on the group ${\mathcal
J}({\bf x})$; this monodromy action is described by the Picard--Lefschetz
formulae, see \S \ 2. Lemma 6 follows from these formulae, from Lemma 5, and
from the fact that the group ${\mathcal J}({\bf x})$ coincides with the linear
hull of the orbit of any vanishing cycle under the action of the big monodromy
group, see [Gab], [E].
\medskip

The set $L \cap \Sigma(F)$ consists of several points of two kinds: the points
of transversal intersection of $L$ and $W_F$ and points $x \not \in W_F$ such
that $S(x)$ is tangent to $W_F$. \medskip

{\bf Lemma 7.} {\it a) Close to a generic point $y$ of the submanifold $W_F
\subset {\C}^n \subset T$ (i.e. to a point at which the generating lines of the
cone $S(y)$ are transversal to $W_F$) the variety $\Sigma_T$ is smooth and has
simple tangency with ${\C}^n$ along $W_F$. In particular, the intersection of
$\Sigma_T$ with any 2-plane $\Lambda$ transversal to $W_F$ coincides close to
the points of $\Lambda \cap W_F$ with a smooth curve having simple tangency
with the line $\Lambda \cap {\C}^n \equiv L$; \smallskip

b) if $F$ is $S$-generic, then close to a generic point of the variety
$(\Sigma(F) \setminus W_F) \subset {\C}^n \subset T$ the variety $\Sigma_T$ is
smooth and intersects ${\C}^n$ transversally along $(\Sigma(F) \setminus
W_F)$.}
\medskip

The proof is immediate. \medskip

Thus the cardinality of $L' \cap \Sigma_T$ is equal to the cardinality of $L
\cap \Sigma(F)$ plus $\; \hbox{deg} \, F$: to any point of $L \cap (\Sigma(F)
\setminus W_F)$ there corresponds one close point of $L' \cap \Sigma_T$, while
to any point of $L \cap W_F$ there correspond two such points; see Fig. 3a.

\begin{figure}
\unitlength=0.80mm \special{em:linewidth 0.4pt} \linethickness{0.4pt}
\begin{picture}(147.00,75.00)
\put(10.00,42.00){\line(1,0){60.00}}
\put(33.00,41.00){\rule{2.00\unitlength}{2.00\unitlength}}
\put(10.00,46.00){\line(6,-1){60.00}} \put(18.00,51.00){\oval(10.00,18.00)[b]}
\put(18.00,42.00){\circle*{2.00}} \put(52.00,33.50){\oval(10.00,17.00)[t]}
\put(52.00,42.00){\circle*{2.00}} \put(58.00,55.00){\line(2,-5){11.67}}
\put(63.00,42.00){\circle{2.00}} \put(45.00,55.00){\line(-1,-6){4.67}}
\put(43.00,42.00){\circle{2.00}} \put(69.00,44.00){\makebox(0,0)[cc]{$L$}}
\put(69.00,39.00){\makebox(0,0)[cc]{$L'$}}
\put(34.00,39.00){\makebox(0,0)[cc]{${\bf x}$}}
\put(34.00,63.00){\makebox(0,0)[cc]{\large $\Lambda$}}
\put(32.00,4.00){\makebox(0,0)[cc]{\large $a$}}
\put(129.00,61.00){\rule{2.00\unitlength}{2.00\unitlength}}
\put(130.00,58.00){\makebox(0,0)[cc]{${\bf x}$}}
\put(140.00,62.00){\circle*{2.00}} \put(145.00,72.00){\circle{2.00}}
\put(145.00,52.00){\circle{2.00}} \put(125.00,50.00){\circle{2.00}}
\put(125.00,74.00){\circle{2.00}} \put(114.00,62.00){\circle*{2.00}}
\put(107.00,62.00){\circle*{2.00}}
\put(129.00,24.00){\rule{2.00\unitlength}{2.00\unitlength}}
\put(131.00,21.00){\makebox(0,0)[cc]{${\bf x}$}}
\put(146.00,35.00){\circle{2.00}} \put(144.00,14.00){\circle{2.00}}
\put(126.00,12.00){\circle{2.00}} \put(123.00,37.00){\circle{2.00}}
\put(130.00,25.00){\line(-1,-3){4.00}} \put(130.00,25.00){\line(-3,5){6.67}}
\put(131.00,26.00){\line(5,3){14.33}} \put(131.00,24.00){\line(4,-3){12.33}}
\put(140.67,27.00){\circle*{1.33}} \put(131.00,25.67){\line(1,0){6.00}}
\put(137.00,25.67){\line(3,1){3.67}} \put(131.00,24.67){\line(1,0){6.00}}
\put(137.00,24.67){\line(5,-3){4.67}} \put(141.67,22.00){\circle*{1.33}}
\put(129.00,25.67){\line(-1,0){10.00}} \put(119.00,25.67){\line(-3,2){4.00}}
\put(115.00,28.33){\circle*{1.33}} \put(129.00,24.67){\line(-1,0){9.67}}
\put(119.33,24.67){\line(-3,-1){5.00}} \put(114.00,22.67){\circle*{1.33}}
\put(129.00,24.67){\line(-3,-1){20.33}} \put(108.67,18.00){\line(0,1){10.33}}
\put(108.67,29.00){\circle*{1.33}} \put(130.00,24.00){\line(-3,-1){22.67}}
\put(107.33,16.33){\line(-1,6){1.33}} \put(106.33,24.00){\circle*{1.33}}
\put(95.00,24.00){\makebox(0,0)[cc]{$L'$}}
\put(95.00,62.00){\makebox(0,0)[cc]{$L$}}
\put(115.00,4.00){\makebox(0,0)[cc]{\large $b$}}
\end{picture}

\centerline{Fig. 3. Lines $L$ and $L'$ and discriminant points in them}
\end{figure}

Since the point ${\bf x}$ lies in the hyperbolicity domain, all points of $L
\cap W_F$ are real. For any such point $y$ belonging to the $k$-th component of
$M_F$, let $y_+ \in {\R}^n \setminus \Sigma(F)$ be a close point in the $k$-th
zone. For such a point $y_+$, the class $a(y_+)$ was defined before Theorem 9.

Let us agree to choose the distinguished system of paths in $L'$ in such a way
that the paths connecting ${\bf x}$ with any two points of $L' \cap \Sigma_T$
arising from the same point $y$ of $L \cap W_F$ go together up to a small
common neighborhood of these two points and are close to the real segment in
$L$ connecting ${\bf x}$ and $y$, while the paths in $L'$ connecting ${\bf x}$
with any other points of $L' \cap \Sigma_T$ do not touch this small
neighborhood; see Fig. 3b. \medskip

{\bf Definition 12.} A point of $L' \cap \Sigma_T$ is of the first kind
(respectively, of the second kind) if it arises from a close point of $W_F$
(respectively, of $\Sigma(F) \setminus W_F$) in $L$ after the move $L \to L'$.
A cycle in ${\mathcal J}({\bf x})$ vanishing  over a path of our distinguished
system in $L'$ that connects ${\bf x}$ with a point $y \in \Sigma_T$ is called
a cycle of the first kind (respectively, of the second kind) if this point $y$
is of the first (respectively, the second) kind. \medskip

In Fig. 3b the points of $L \cap W_F$ and the points of the first kind in $L'$
are shown by small black circles, while the points of $L \cap (\Sigma(F)
\setminus W_F)$ and the points of the second kind in $L'$ are shown by white
circles.
\medskip

{\bf Lemma 8.} {\it a) Two cycles of the first kind in ${\mathcal H}({\bf x})$,
vanishing over two distinguished paths connecting ${\bf x}$ with two points of
$L' \cap \Sigma_T$ arising from the same close point $y$ of $L \cap W_F$,
coincide (maybe up to sign); \smallskip

b) this cycle coincides (maybe up to sign) with the cycle $a(y_+)$ transported
from the point $y_+$ to ${\bf x}$ along the path described in the definition of
the reduced Arnold class. In particular, the map $\Pi_*$ sends the homology
class of any such cycle into a generator of the corresponding group (\ref{pr1})
or (\ref{pr2}); \smallskip

c) the monodromy action in the group ${\mathcal H}({\bf x})$, defined by any
simple loop in $L \setminus \Sigma(F)$ going around some point of $L \cap W_F$,
is trivial; \smallskip

d) any cycle in ${\mathcal H}({\bf x})$ vanishing over a path in $L \setminus
\Sigma(F)$ connecting ${\bf x}$ with a point of $\Sigma(F) \setminus W_F$
belongs to the subspace ${\mathcal M}({\bf x})$. In particular, the same is
true for any cycle of the second kind defined by a path of our distinguished
system in $L' \setminus \Sigma_T$ connecting ${\bf x}$ with a point (of the
second kind) of $\Sigma_T$.} \medskip

{\it Proof.} Consider the space of complex lines through ${\bf x}$ transversal
to $\Sigma_T$ in the plane $\Lambda$. Obviously this space is a projective line
with several points removed, one of which is the point $\{L\}$. Consider a
small loop in this space, which starts and finishes at the point $\{L'\}$ and
goes once around the point $\{L\}$. This loop takes one of the two
distinguished paths from statement a) of the lemma into the other, thus this
statement follows.

Statement c) is a direct consequence of a). Indeed, the loop considered there
is homotopic in $\Lambda \setminus \Sigma_T$ to a loop in $L' \setminus
\Sigma_T$ which turns around two discriminant points defining the same
vanishing cycle, thus its monodromy action is equal to the square of the
reflection in the hyperplane orthogonal to this vanishing cycle.

Statement b) follows from Lemma 4 and the local shape of the pair $(W_F$,
${\mathcal S}(\lambda))$ where ${\mathcal S}(\lambda)$ is the variety of zeros
of the polynomial (\ref{ttt}) defined by the discriminant point $\lambda$ of
the first kind.  The way in which the pairs of distinguished paths  connecting
$x$ with different pairs of points of the first kind miss one another is not
important, because by the proof of Theorem 3 all the cycles of the first kind
that vanish over the paths going from ${\bf x}$ to the points arising from
different points of $L\cap W_F$ on  the same side of ${\bf x}$ in $\; \hbox{Re}
\, L$ are pairwise orthogonal.

Statement d) of the lemma follows immediately from the constructions.
\medskip

Thus, the vanishing cycles of the first (respectively, second) kind are exactly
those that are sent by the map $\Pi_*$ into a generator of the group
(\ref{pr1}) or (\ref{pr2}) (respectively, into a zero class). \medskip

{\bf Lemma 9.} {\it Any vanishing cycle of the first kind in ${\mathcal J}({\bf
x})$ can be transformed into any other by a sequence of reflections in the
hyperplanes orthogonal to cycles of the second kind and to this cycle itself.}
\medskip

By the Picard--Lefschetz formula, this lemma follows from the next one.
\medskip

{\bf Lemma 9$'$.} {\it There exists a distinguished system of paths in $L'
\setminus \Sigma_T$ connecting ${\bf x}$ with all points of $L' \cap \Sigma_T$,
such that all vanishing cycles of the first kind defined by this system are
equal to each other.} \medskip

{\it Proof.} (This proof simulates that of the well-known fact that the
fundamental group of the complement of a smooth irreducible algebraic
hypersurface in ${\C}^n, \, n\ge 2$ is isomorphic to ${\ZZ}$.)

Let $y_1$ be any point of $L \cap W_F$. Let us fix an arbitrary path $\gamma_1$
in $L \setminus \Sigma_T$ connecting ${\bf x}$ with $y_1$. Denote by $A^n$ the
space of complex lines in ${\C}^n$, and by $\hbox{Reg}\, (\Sigma(F))$ the
subset of $A^n$ consisting of lines transversal to $\Sigma(F)$. Consider a path
$\chi_1:[0,1] \to A^n$ such that $\chi_1(0) = L$, $\chi_1([0,1)) \subset
\hbox{Reg}\,(\Sigma(F))$, the last point $\chi_1(1)$ is a line transversal to
$\Sigma(F)$ everywhere except for one point of simple tangency with $W_F$, and
one of the two points of $\chi_1(\tau) \cap W_F, \; \tau =1-\varepsilon$ that
coalesce at this tangency point is obtained from the point $y_1$ of the similar
set corresponding to the value $\tau = 0$ during the deformation of the set
$\chi_1(\tau) \cap W_F, \; \tau \in [0,1-\varepsilon]$.

Consider the continuous deformation $\gamma_1[\tau], \,\tau \in [0,1]$, of the
path $\gamma_1$ such that $\gamma_1[0] = \gamma_1$, $\gamma_1[\tau] \subset
\chi_1(\tau)$, and for any $\tau$ the path $\gamma_1[\tau]$ connects in
$\chi_1(\tau) \setminus \Sigma(F)$ a point of $\chi_1(\tau) \cap W_F$ with some
distinguished point ${\bf x}(\tau) \in \gamma_1[\tau] \setminus \Sigma(F)$,
${\bf x}(0)={\bf x}$. At  almost the final instant $\tau=1- \varepsilon$, the
endpoint $\gamma_1[1-\varepsilon](1)$ of the path $\gamma_1[1-\varepsilon]$
lies very close to some other point of $\chi_1(1-\varepsilon)\cap W_F$ (with
which it coalesces at the instant $\tau = 1$). Connect this new point with
${\bf x}(1-\varepsilon)$ by a path $\gamma_2[1-\varepsilon]$ in
$\chi_1(1-\varepsilon)\setminus \Sigma(F)$ that goes very close to
$\gamma_1[1-\varepsilon]$ but does not intersect it except for the initial
point. Then construct a continuous family of paths $\gamma_2[\tau] \subset
\chi_1(\tau), \,\tau \in [0,1-\varepsilon]$, such that for any $\tau$ the
corresponding path $\gamma_2[\tau]$ connects a point of $\chi_1[\tau] \cap W_F$
with ${\bf x}(\tau)$ and does not intersect other points of $\chi_1(\tau) \cap
\Sigma(F)$ or of the path $\gamma_1[\tau]$. At the instant $\tau=0$ we get a
path $\gamma_2 \equiv \gamma_2[0] \subset L$ connecting ${\bf x}$ with some
point $y_2$ of $W_F$.

Then consider a new path $\chi_2:[0,1] \to A^n$, $\chi_2([0,1)) \subset
\hbox{Reg} \,(\Sigma(F))$, connecting $L$ with some new simple tangent line to
$W_F$ and having no extra nontransversalities with $\Sigma(F)$, in such a way
that at the last instant $\tau=1$ one of the two points of $\chi_2(\tau) \cap
W_F$ that  coalesce at the tangency point is obtained by deformation along our
path $\chi_2$ from one of the points $y_1$ or $y_2$, and the other two points
of these two pairs do not coincide. Arguing as before, we construct a third
path in $L \setminus \Sigma(F)$, connecting ${\bf x}$ with some third point of
$L \cap W_F$, and so on.

After the $(d-1)$-th step we get a system of $d$ nonintersecting paths in $L
\setminus \Sigma(F)$, connecting ${\bf x}$ with all points of $L \cap W_F$.
Complete this family to any distinguished collection of paths connecting ${\bf
x}$ with all points of $L\cap \Sigma(F)$. For the close perturbation $L'
\subset T$ of $L$, take a close distinguished system of paths in $L'$,
connecting the point ${\bf x}$ with all points of $L' \cap \Sigma_T$ in such a
way that to any path in $L$ connecting ${\bf x}$ with $W_F$ there correspond
two paths connecting ${\bf x}$ with two close points of the first kind. This
system of paths is the desired one. For instance, the cycles vanishing along
the (perturbed) paths $\gamma_1$ and $\gamma_2$ define the same vanishing
homology class in ${\mathcal J}({\bf x})$: indeed, a similar assertion for the
cycles in the group ${\mathcal H}({\bf x}(1-\varepsilon)) \equiv  H_{n-1}(W_F
\setminus S({\bf x}(1-\varepsilon)), {\ZZ})$ or $ H_{n-1}(W_F \setminus S({\bf
x}(1-\varepsilon)), @({\bf x}(1-\varepsilon)))$ is proved just as the statement
a) of Lemma 8, and for other values of $\tau \in [0,1-\varepsilon]$ it follows
by continuity. Lemmas $9'$ and 9 are thus proved.
\medskip

Now we are ready to prove statement b) of Theorem 7. Indeed, by Lemma 6 the
group ${\mathcal J}({\bf x})$ is generated by the vanishing cycles of the first
and second kind. By Lemma 9 and the Picard--Lefschetz formula, all vanishing
cycles of the first kind lie in the linear span of an arbitrary one of them
(for which we can take the class obtained by the Gauss--Manin connection from
$a(x)$, $x$ from the ${\bf k}$-th zone, $1 \le k \le [d/2]$, see statement b)
of Lemma 8) and the vanishing cycles of the second kind (which lie in
${\mathcal M}({\bf x})$, see statement d) of Lemma 8). \smallskip

Statement c) of Theorem 7 follows immediately from statement c) of Lemma 8, and
statement d) follows from the fact that the variety $\Sigma(F) \setminus W_F$
is irreducible. \medskip

{\it Proof of Theorem 9a).} We can assume that the points $y_1$ and $y_2$,
whose classes $a(y_1)$ and $a(y_2)$ we want to transfer to each other, lie very
close to the ``interior'' (i.e. closest to the hyperbolicity domain) components
of $M_F$ bounding corresponding zones. For such $y_i$ the class $a(y_i)$ is
realized by a cycle generating the group $ H_{n-1}(W_F \cap B \setminus
S(y_i))$ or $ H_{n-1}(W_F \cap B \setminus S(y_i), @(y_i))$, where $B$ is a
small neighbourhood of $y_i$; see Lemma 4. Thus, for the desired path
connecting $y_1$ and $y_2$ we can take the path that goes very close to the set
of generic points of $W_F$ (i.e. of such points $y$ close to which all the
generating lines of the cones $S(y)$ are transversal to $W_F$ and hence the
pairs  $(W_F, S(y))$ have locally the same topological structure). \smallskip

Statement b) of Theorem 9 follows from Theorem 3 and the connectedness of
Dynkin diagrams of isolated singularities of complete intersections.
\medskip

{\it Proof of the statement g) of Theorem 7.}

First of all, this statement is true in the case when $M_F$ is an ellipsoid
with different eigenvalues. Indeed, by Theorem 3 in this case $\tilde A(x)$ is
a vanishing cycle, and the assertion follows from the connectedness of the
Dynkin diagram and the fact that the group ${\mathcal M}(x)$ is nontrivial for
such $F$, see e.g. [E].

For arbitrary $d$, consider the model (not $S$-generic) hyperbolic surface
$M_{F}$ consisting of $[d/2]$ ellipsoids $ \alpha_1 x_1^2 + \cdots + \alpha_n
x_n^2 = j, \quad j=1, 1+\varepsilon, \ldots, 1+([d/2]-1)\varepsilon, $ where
all $\alpha_i$ are positive and distinct, plus, if $d$ is odd, a distant
hyperplane. The class $\tilde A(x)$ for $x$ from the $k$-th zone, $1\le k \le
[d/2]$, is then equal to the sum of $k$ vanishing cycles, each of which lies in
the complexification of its own ellipsoid; see the proof of Theorem 3. By the
previous special case of a single ellipsoid, in each of these $k$ complexified
ellipsoids ${\mathcal E}_i$ there is a compact cycle $\Gamma$ defining an
element of the group $ H_{n-1}({\mathcal E}_i \setminus S(x))$ if $n$ is even,
or in $ H_{n-1}({\mathcal E}_i \setminus S(x), @(x))$ if $n$ is odd, such that
$\langle \tilde A(x), \Gamma \rangle \ne 0$ and the map $\Pi_*$ sends the
homology class of $\Gamma$ into the zero homology class.

Consider a perturbation of our model hyperbolic polynomial $F$ which replaces
it by a $S$-generic one and is so weak that it does not change the topology of
the variety $W_{F} \cup S(x)$ inside a sufficiently large disc, in which the
cycles $\Gamma$  and $\tilde A(x)$ lie. The cycle $\tilde \Gamma$ close to
$\Gamma$ in the moved manifold $W_F$ satisfies all the above conditions, and
statement g) of Theorem 7 is proved for {\em some} $S$-generic hyperbolic
polynomial. For an arbitrary such polynomial this statement follows from the
fact that all the generic surgeries separating different path-components of the
space of all strictly hyperbolic $S$-generic surfaces of given degree in $\R^n$
(these surgeries correspond to the smooth hyperbolic surfaces in $\R P^n$
simple tangent to the non-proper plane) preserve the homology classes $\tilde
A(x)$ (provided that the corresponding point $x$ and the distinguished point in
the hyperbolicity domain do not change in this surgery). (In formal terms, this
preservation means that these homology classes corresponding to the polynomials
before and after the surgery are transposed into one another by the natural
connection over any {\em short} connecting them path in the space of all {\em
complex} $S$-generic polynomials.)

\subsection{Proof of Theorem 8}

Let $c$ be a point of simple tangency of a cone $S(x_0)$ and $W_F$ such that
$P(c) \ne 0$. Let ${\Upsilon}$ be an affine complex line through $x_0$ in
$\C^n$, transversal to the common tangent hyperplane of $S(x_0)$ and $W_F$ at
$c$; let $\xi$ be an affine coordinate on it with the origin at $x_0$. Consider
the one-parametric family of surfaces $S(x(\xi)),$ $x(\xi) \in {\Upsilon}$. The
elements $S(x(\xi))$ of this family with $\xi$ from a small {\it punctured}
neighborhood of the origin are transversal to $W_F$ in a small disc $B$ centred
at $c$, and the vanishing element $\gamma(\xi)$ of the group $H_{n-1}(B \cap
W_F \setminus S(x(\xi)))$ (if $n$ is even) or $H_{n-1}(B \cap W_F \setminus
S(x(\xi)), @(x(\xi)))$ (if $n$ is odd) is well defined (up to sign) by this
family.  By the Picard--Lefschetz formulae of \S \ 2, in both cases the
rotation of $\xi$ around $0$ sends $\gamma(\xi)$ to $-\gamma(\xi)$.

Define the function $\Xi(\xi), \, $ $\xi \in \C$, as the integral of the form
(\ref{polpot}) with $x=x(\xi) \in \Upsilon$ along the cycle $\gamma(\xi)$. It
is sufficient to prove that there are arbitrarily high derivatives of this
function not equal identically to 0.  This follows from the next lemma.
\medskip

{\bf Lemma 11.} {\it The function $\Xi(\xi)$ is represented by a power series
of the variable $\sqrt{\xi}$, whose leading (of smallest degree) term with
non-zero coefficient has degree $1$.}

(Of course, all even powers of this series vanish.)
\medskip

{\it Proof.} Using the Leray residue theorem, we can replace the integral
(\ref{polpot}) along the cycle $\gamma(\xi)$ by the integral of the form
$G(x(\xi)-z)P(z)/F(z)dz_1 \wedge \ldots \wedge dz_n$ along the Leray tube $t
\gamma(\xi) \in H_n(B\setminus (W_F\cup S(x(\xi))))$ or $\in H_n(B\setminus
(W_F\cup S(x(\xi))), @(x(\xi)))$. Close to $c$ the holomorphic function $\C^n
\to \C$ is defined, which assigns to any point the coordinate $\xi$ of the
origin $x(\xi)$ of the cone $S(x(\xi))$ containing it. Choose this function for
the last local coordinate $w_n$ at $c$; by the Morse lemma we can choose the
remaining coordinates $w_1, \ldots, w_{n-1}$ in such a way that $W_F$ is
locally given by $w_n = w_1^2 + \cdots + w_{n-1}^2.$ In these coordinates our
differential form becomes
\begin{equation}
\label{integ} (w_n-\xi)^{-(n-2)/2}(w_n - w_1^2 - \cdots - w_{n-1}^2)^{-1}
I(w_1, \ldots, w_n)dw_1 \wedge \ldots \wedge dw_n,
\end{equation}
where the function $I$ does not vanish at $c$. Let $I= I_0 + I_1 + \cdots $ be
the expansion of $I$ into the sum of quasihomogeneous polynomials of degrees
$0, 1 , \ldots$ respectively with respect to the weights $\deg w_1 = \cdots =
\deg w_{n-1}=1, \deg w_n=2.$ Using the corresponding group of quasihomogeneous
dilations $(w_1, \ldots, w_{n-1}, w_n) \to (\tau w_1, \ldots, \tau w_{n-1},
\tau^2w_n)$ we see, that the integral along $t \gamma(\xi)$ of the form similar
to (\ref{integ}), in which $I_m$ is substituted instead of $I$, is a
homogeneous function in $\xi$ of degree $(m+1)/2$. It is easy to calculate that
this function corresponding to the constant polynomial $I_0 \ne 0$ is not the
identical zero function; this proves our Lemma.


\begin{thebibliography}{999}

\bibitem [{\bf A 82}] Z V.~I.~Arnold,
{\it On the Newtonian potential of hyperbolic layers,} Trudy Tbilisskogo
Universiteta, Ser. Math./Mekh./Astron., 1982, V. 232--233, \# 13--14, 23--29.
English transl. in: Selecta Math. Soviet., 4:2, 1985, 103--106.

\bibitem [{\bf A 83}] Z V.~I.~Arnold, {\it Magnetic field analogues
of the theorems of Newton and Ivory.} Uspekhi Mat. Nauk 1983, 38:5, 145--146.

\bibitem [{\bf A 87}] Z V.~I.~Arnold, {\it Kepler$'$s second
law and the topology of Abelian integrals}, Kvant, No. 12, 17--21 (Russian);
{\it A topological proof of transcendence of Abelian integrals in Newton$'$s
Principia}, Quant, 1987, No.12, 1--15; {\it The 300-th anniversary of
mathematical natural philosophy and celestial mechanics}, Priroda, 1987, No.8,
5--15 (in Russian).

\bibitem [{\bf ABG}] Z M.~F.~Atiyah, R.~Bott, L.~G\aa rding,
{\it Lacunas for hyperbolic differential operators with constant coefficients},
Acta Math., 124 (1970), 109--189 and 131 (1973),  145--206.

\bibitem [{\bf AV}] Z V.~I.~Arnold and V.~A.~Vassiliev, {\it Newton$'$s
Principia read 300 years later}, Notices Amer. Math. Soc., 36:9 (1989),
1148--1154.

\bibitem [{\bf AVG}] Z V.~I.~Arnold, A.~N.~Varchenko, S.~M.~Gusein-Zade,
{\it Singularities of differentiable maps}, V 1, 2; ``Nauka'', Moscow, 1982 and
1984; Engl. transl.:  Birkh\"auser, Basel, 1985 and 1988.

\bibitem [{\bf AVGL}] Z V.~I.~Arnold, V.~A.~Vassiliev, V.~V.~Goryunov,
O.~V.~Lyashko, {\it Singularities, 1 and 2}. Dynamical systems, 6 and 39,
VINITI, Moscow, 1988 and 1989; English transl.: Encyclopaedia Math. Sci., V. 6
and 39, Springer-Verlag, Berlin and New York, 1993.

\bibitem [{\bf E}] Z  W.~Ebeling, {\it The monodromy groups of
isolated singularities of complete intersections.} Lect. Notes Math., vol.1293,
Springer, Berlin a.o., 1987.

\bibitem [{\bf G 84}] Z A.~B.~Givental, {\it Polynomiality of the
electrostatic potentials,} Uspekhi Mat. Nauk, 39:5 (1984), 253--254 (in
Russian).

\bibitem [{\bf G 88}] Z A.~B.~Givental, {\it Twisted Picard--Lefschetz
formulae.} Funkts. Anal. i Prilozh., 22:1 (1988), 12--22; Engl. translation in
Functional Anal. Appl., 22:1, 10--18 (1988).

\bibitem [{\bf Gab}]{Gab} A.~M.~Gabrielov, {\it
Bifurcations, Dynkin diagrams and modality of isolated singularities,} Funkts.
Anal. i Prilozh., 8:2 (1974), 7--12; Engl. transl. in Funct. Anal. Appl., 8
(1974), 94--98.

\bibitem [{\bf GH}] Z G.~M.~Greuel, H.~A.~Hamm, {\it
Invarianten quasihomogener voll\-st\"an\-diger Durchschnitte}, Invent. Math.,
1978, 49:1, 67--86.

\bibitem [{\bf GrH}] Z Ph.~Griffiths, J.~Harris, {\it
Principles of algebraic geometry}, John Wiley \& Sons, New York a.o., 1978.

\bibitem[{\bf GM}]{GM}
M.~Goresky and R.~MacPherson, {\it Stratified Morse theory}, Springer-Verlag,
Berlin and New York, 1986.

\bibitem [{\bf GR}] Z H.~Grauert and R.~Remmert,
{\it Komplexe R\"aume}, Math.  Ann., 136:2 (1958), 245--318.

\bibitem [{\bf GZ}] Z S.~M.~Gusein--Zade, {\it Monodromy groups of
isolated singularities of hypersurfaces}, Uspekhi Mat. Nauk 32 (1977), no.2,
23--65, Engl. transl. in Russian Math. Surveys 32, No.2 (1977), 23--69.

\bibitem [{\bf H}] Z H.~Hamm, {\it Locale topologische Eigenschaften
komplexer R\"aume,} Math. Ann., 191, 1971, 235--252.

\bibitem [{\bf I}] Z J.~Ivory, {\it On the attraction of homogeneous
ellipsoids}, Philos. Trans., 99 (1809), 345--372.

\bibitem [{\bf M}] Z J.~Milnor,
{\it Singular points of complex hypersurfaces}, Princeton Univ. Press,
Princeton, NJ, and Univ. of Tokyo Press, Tokyo, 1968.

\bibitem [{\bf MO}] Z J.~Milnor, P.~Orlik, {\it Isolated singularities,
defined by weighted isolated polynomials}, Topology, 1970, 9:2, 385--393.

\bibitem [{\bf Newton}] Z I.~Newton, {\it Philosophiae Naturalis Principia
Mathematica}, London, 1687.

\bibitem [{\bf N}] Z W.~Nuij, {\it A note on hyperbolic polynomials},
Math. Scand. 23 (1968), 69--72.

\bibitem [{\bf P}] Z I.~G.~Petrovsky, {\it On the diffusion of waves
and the lacunas for hyperbolic equations}, Matem. Sbornik, 1945, 17(59),
289--370.

\bibitem [{\bf Ph 65}] Z F.~Pham, {\it Formules de Picard--Lefschetz
g\'en\`eralis\'ees et ramification des in\-t\'e\-gra\-les}, Bull. Soc. Math.
France, 93 (1965), 333--367.

\bibitem [{\bf Ph 67}] Z F.~Pham, {\it Introduction \`a l'\'etude topologique
des singularit\'es de Landau}, Gauthier-Villars, Paris, 1967.

\bibitem [{\bf VSh}] Z A.~D.~Vainshtein and B.~Z.~Shapiro,
{\it Multidimensional analogues of the Newton and Ivory theorems}, Funkts.
Anal. i Prilozh., 19:1 (1985), 20--24; Engl. translation in Functional Anal.
Appl., 19:1 (1985), 17--20.

\bibitem [{\bf V 86}] Z V.~A.~Vassiliev, {\it Sharpness and
the local Petrovskii condition for hyperbolic equations with constant
coefficients,} Izv. Akad. Nauk SSSR Ser. Mat. 50 (1986), 242--283; Engl.
transl. in Mat. USSR Izv. 28 (1987), 233--273.

\bibitem [{\bf V 94}] Z V.~A.~Vassiliev, {\it
Ramified Integrals, Singularities and Lacunas}, Kluwer, 1994.

\bibitem [{\bf Z}] Z O.~Zariski, {\it On the Poincar\'e group of a
projective hypersurface}, Ann. Math. 38 (1937), 131--141.

\end{thebibliography}
\end{document}